# EXCLUSION PROCESSES IN HIGHER DIMENSIONS: STATIONARY MEASURES AND CONVERGENCE

By M. Bramson[1] and T. M. Liggett[2]

*University of Minnesota and University of California, Los Angeles*

There has been significant progress recently in our understanding of the stationary measures of the exclusion process on $Z$. The corresponding situation in higher dimensions remains largely a mystery. In this paper we give necessary and sufficient conditions for a product measure to be stationary for the exclusion process on an arbitrary set, and apply this result to find examples on $Z^d$ and on homogeneous trees in which product measures are stationary even when they are neither homogeneous nor reversible. We then begin the task of narrowing down the possibilities for existence of other stationary measures for the process on $Z^d$. In particular, we study stationary measures that are invariant under translations in all directions orthogonal to a fixed nonzero vector. We then prove a number of convergence results as $t \to \infty$ for the measure of the exclusion process. Under appropriate initial conditions, we show convergence of such measures to the above stationary measures. We also employ hydrodynamics to provide further examples of convergence.

**1. Introduction.** In this paper we consider the stationary measures and convergence for the exclusion process $\eta_t$ on the countable set $S$ with (stochastic) transition probabilities $p(x, y)$. This is the Markov process on $\{0, 1\}^S$ in which a particle at $x \in S$ attempts to move to $y \in S$ at rate $p(x, y)$. The move takes place if $y$ is vacant and does not take place if $y$ is occupied. The exclusion process has been very heavily studied since its introduction by Spitzer [18] in 1970. For an account of many known results about it, see the books by Liggett [13, 14], and Kipnis and Landim [8], and the references therein.

In spite of the tremendous success that the study of the exclusion process has seen, relatively little was known about the existence and structure

Received March 2004; revised October 2004.
[1]Supported in part by NSF Grant DMS-02-26245.
[2]Supported in part by NSF Grant DMS-03-01795.

of nonexchangeable and nonreversible stationary measures for the translation invariant exclusion process on $Z^d$ until 2001. (The only truly general result that was known is that all extremal translation invariant stationary measures are of product form—see Theorem 3.9 of Chapter VIII in [13].) This is still the case if $d > 1$. In one dimension, however, the papers by Ferrari, Lebowitz and Speer [5], Bramson and Mountford [2] and Bramson, Liggett and Mountford [1] made significant progress on this problem. The present paper is a first attempt to develop a corresponding theory in higher dimensions.

We begin by summarizing some of the one-dimensional results, since they provide a model for what we would like to prove in higher dimensions. The main assumptions are $S = Z$, $p(x,y) = p(y-x)$, $\sum_x |x| p(x) < \infty$, and the process is irreducible. [Throughout this paper, irreducibility will mean that, for each $x, y \in S$, there is positive probability that the Markov chain with transition probabilities $p(\cdot, \cdot)$ goes from $x$ to $y$ in some number of steps.] First, two definitions: A probability measure $\mu$ on $\{0,1\}^Z$ is said to be a *blocking measure* if it concentrates on configurations $\eta$ satisfying

$$\sum_{x<0} \eta(x) < \infty \quad \text{and} \quad \sum_{x>0} [1 - \eta(x)] < \infty,$$

and it is said to be a *profile measure* if

$$\lim_{x \to -\infty} \mu\{\eta : \eta(x) = 1\} = 0 \quad \text{and} \quad \lim_{x \to +\infty} \mu\{\eta : \eta(x) = 1\} = 1.$$

Every blocking measure is a profile measure, but not conversely. As usual, if $0 \leq \rho \leq 1$, $\nu_\rho$ will denote the homogeneous product measure with density $\rho$.

Without loss of generality, we may assume $\sum_x x p(x) \geq 0$. We then have the following results:

1. Either the extremal stationary measures consist exactly of

(i) $$\{\nu_\rho, 0 \leq \rho \leq 1\}$$

or of

(ii) $$\{\nu_\rho, 0 \leq \rho \leq 1\} \cup \{\mu_n, n \in Z\},$$

where $\mu_0$ is a profile measure, and $\mu_n$ is the shift of $\mu_0$ by $n$ units.
2. If $\sum_x x p(x) = 0$, then case (i) occurs.
3. If $\sum_x x p(x) > 0$ and $p(\cdot)$ has finite range, then case (ii) occurs and $\mu_0$ is a blocking measure.
4. If $\sum_x x p(x) > 0$; $p(x)$ and $p(-x)$ are decreasing for $x = 1, 2, \ldots$; $p(x) \geq p(-x)$ for all $x \geq 1$; and $\sum_{x<0} x^2 p(x) < \infty$, then case (ii) occurs and $\mu_0$ is a blocking measure.
5. If $\sum_{x<0} x^2 p(x) = \infty$, there exists no stationary blocking measure.

Statements 1, 4 and 5 were proved in [1], statement 2 appears on page 391 in [13], and statement 3 is due to [2]. An important open problem that remains in one dimension is to determine whether nonblocking stationary profile measures ever exist.

While the above description is the starting point for our study of stationary measures in higher dimensions, one difference between one and higher dimensions can be seen immediately in case (ii) above. In one dimension, the extremal spatially inhomogeneous stationary measures are naturally indexed by a discrete parameter. In higher dimensions, the natural parametrization is continuous, rather than discrete. Here is a simple example in two dimensions that illustrates this point. Suppose $p(e_1) = p(e_2) = p(-e_2) = 1/3$, where $e_1 = (1,0)$ and $e_2 = (0,1)$. Then, for each $s \in R$, an extremal stationary measure $\mu_s$ is obtained in the following way. Letting $x = (x^{(1)}, x^{(2)})$ and $s = n + \sigma$ for an integer $n$ and $0 \le \sigma < 1$, $\mu_s$ is the measure for which $\eta(x) \equiv 1$ if $x^{(1)} > n$, $\eta(x) \equiv 0$ if $x^{(1)} < n$, and $\eta(x)$ is i.i.d. with density $\sigma$ if $x^{(1)} = n$.

We turn now to a description of the results of the current paper. Precise statements and definitions appear later in the relevant sections. For the moment, $S$ is general. The first problem is to determine exactly when product measures are stationary for the process. Suppose that $\alpha(x)$ satisfies $0 < \alpha(x) < 1$ for all $x \in S$, set

$$\pi(x) = \alpha(x)/(1 - \alpha(x)),$$

and let $\nu_\alpha$ be the product measure on $\{0,1\}^S$ with marginals $\nu_\alpha\{\eta : \eta(x) = 1\} = \alpha(x)$. It is not hard to check that $\nu_\alpha$ is reversible for $\eta_t$ if and only if $\pi(x)p(x,y) = \pi(y)p(y,x)$ for all $x, y \in S$—see page 34 in [6], for example.

Here we are interested in the issue of stationarity of $\nu_\alpha$, rather than reversibility. Theorem 2.1 on page 380 in [13] says that $\nu_\alpha$ is stationary for $\eta_t$ if either (a) $p(x,y)$ is doubly stochastic and $\alpha(x)$ is constant, or (b) $\pi(x)p(x,y) = \pi(y)p(y,x)$ for all $x, y \in S$. As mentioned above, in the latter case, $\nu_\alpha$ is in fact reversible for the process. In the former case, $\nu_\alpha$ is reversible if and only if $p(x,y)$ is symmetric. In Section 2 we will show that an appropriate combination of conditions (a) and (b) is necessary and sufficient for $\nu_\alpha$ to be stationary for $\eta_t$.

Section 3 applies this necessary and sufficient condition to the case of a translation invariant system on $S = Z^d$. Under a mild assumption, Theorems 2 and 3 combine to show that $\nu_\alpha$ is stationary if and only if there is a $v \in R^d$ so that $\pi(x) = \pi(0)e^{\langle x,v \rangle}$ and $p(z) = e^{\langle z,v \rangle}p(-z)$ for all $z$ such that $\langle z,v \rangle \ne 0$. This allows us to construct many examples of stationary product measures that are neither homogeneous nor reversible. (Homogeneity of $\nu_\alpha$ corresponds to $\alpha$ being constant.)

In Section 4 we begin by considering several examples in which $p(x,y)$ are the transition probabilities for a random walk on a homogeneous tree $S$.

Again, using the results of Section 2, we find examples of stationary product measures that are neither homogeneous nor reversible. Then, we describe an example on a rooted tree that exhibits a type of phase transition, where the extremality of the stationary product measures among stationary measures depends on the value of a parameter.

In Sections 5 and 6 we return to $S = Z^d$ and study stationary measures that are invariant under translations in all directions orthogonal to a fixed nonzero $v \in Z^d$. We suspect that all extremal stationary measures have this property for some $v \in R^d$, but are quite far from being able to prove this. A measure that has this partial homogeneity property and has density tending to 0 in the $-v$ direction and to 1 in the $+v$ direction will be called a *v-profile measure*. Theorem 4 in Section 5 asserts that there are no stationary $v$-profile measures unless $v$ and the drift vector $\sum_x x p(x)$ have nonnegative inner product. (Under an irreducibility assumption, "nonnegative" can be replaced by "positive.") Other approaches to this and related results under varying assumptions appear in Theorem 5 of Section 5, Corollary 2 of Section 6 and Theorem 10 of Section 8 (the last based on hydrodynamics). This rules out half of the possible directions $v$ in which stationary $v$-profile measures can exist.

Theorem 6 in Section 6 provides a way of showing that, if there exists a continuous family of extremal stationary $v$-profile measures, then there can be no other extremal stationary $v$-profile measures. This is applied in Corollary 1 to examples from Section 3 in which families of stationary product measures are known to exist. Theorem 7 of that section provides some reasonable assumptions under which any extremal stationary $v$-profile measure is a blocking measure, in the sense that it concentrates on configurations $\eta$ such that

$$\sum_{k<0} \eta(kv) < \infty \quad \text{and} \quad \sum_{k>0} [1 - \eta(kv)] < \infty.$$

Theorem 8 of Section 7 provides sufficient conditions for convergence to a mixture of stationary $v$-profile measures when the initial measure is $v$-profile. The main assumption (besides the existence of an appropriate one-parameter family of stationary $v$-profile measures) is that the initial measure have a finite expected number of 1's in the "negative half" of each strip oriented in the $v$ direction (with bounded cross section in the perpendicular directions), and a finite expected number of 0's in the positive half of that strip.

Section 8 presents some applications of hydrodynamical results of Rezakhanlou [16] to obtain explicit convergence results when the initial set of 1's is a wedge in $Z^2$. In this case, the initial distribution is not $v$-profile. Another example with such an initial distribution is discussed in Section 9. It provides an illustration of how the asymptotic behavior of the process can depend on more than just the mean of its transition probabilities.

This paper provides only an initial understanding of the stationary measures and of convergence behavior for the exclusion process in more than one dimension. There are a number of problems in this context that we are not able to solve at this time. These are discussed in Section 10. There, we will also compare our results and speculations with the much easier and well-understood situation for independent particle systems on $Z^d$.

**2. Product stationary measures: the general case.** In this section both $S$ and $p(x,y)$ are general. Our objective is to give a necessary and sufficient condition for a product measure to be stationary for the exclusion process. This is done in the following theorem.

THEOREM 1. *Let $\alpha_1, \alpha_2, \ldots$ be the distinct values of $\alpha(x)$, $\pi_i = \alpha_i/(1-\alpha_i)$, and $C_1, C_2, \ldots$ be the partition of $S$ determined by $C_i = \{x \in S : \alpha(x) = \alpha_i\}$. The following is a necessary and sufficient condition for $\nu_\alpha$ to be stationary for $\eta_t$:*

(a) $\sum_{y \in C_i} p(x,y) = \sum_{y \in C_i} p(y,x)$ *for all $x \in C_i$ and all $i$,*
*and*
(b) $\pi_i p(x,y) = \pi_j p(y,x)$ *for all $x \in C_i, y \in C_j$ and all $i \neq j$.*

PROOF. Following the proof of Theorem 2.1 on page 380 in [13] up to display (2.7), we see that $\nu_\alpha$ is stationary if and only if the right-hand side of (2.7) is zero for every finite $A \subset S$. Dividing by the product $\prod_{u \in A} \alpha(u)$, we see that a necessary and sufficient condition for stationarity is that

$$(1) \quad \sum_{x \in A, y \notin A} \frac{\alpha(y)[1-\alpha(x)]p(y,x) - \alpha(x)[1-\alpha(y)]p(x,y)}{\alpha(x)} = 0$$

for every finite $A \subset S$. Now fix a $u \in S$ and consider (1) with $A = \{u\}$. This becomes

$$(2) \quad [1-\alpha(u)]\sum_y \alpha(y)p(y,u) - \alpha(u)\sum_y [1-\alpha(y)]p(u,y) = 0.$$

(Note that the terms corresponding to $y = u$ in the above sums, which in principle should not be included, cancel out.) Next, fix distinct $v, w \in S$ and consider (1) with $A = \{v, w\}$. This becomes

$$(3) \quad \sum_{y \neq v,w} \frac{\alpha(y)[1-\alpha(v)]p(y,v) - \alpha(v)[1-\alpha(y)]p(v,y)}{\alpha(v)}$$
$$+ \sum_{y \neq v,w} \frac{\alpha(y)[1-\alpha(w)]p(y,w) - \alpha(w)[1-\alpha(y)]p(w,y)}{\alpha(w)} = 0.$$

Subtracting (2) with $u = v$ and (2) with $u = w$ from (3) [after dividing these expressions by $\alpha(v)$ and $\alpha(w)$ resp.], (3) becomes

(4) $\qquad [\alpha(w) - \alpha(v)][\pi(w)p(w,v) - \pi(v)p(v,w)] = 0.$

Therefore, (1) holds for all singleton and doubleton $A$'s if and only if

(5) $\qquad$ (2) holds for all $u$ and (4) holds for all $v, w$.

Next, we will check that (5) implies (1) for all $A$, which gives us the equivalence of (1) and (5). Let $A$ be any finite subset of $S$. By (2), the left-hand side of (1) is

$$-\sum_{x,y \in A} \frac{\alpha(y)[1-\alpha(x)]p(y,x) - \alpha(x)[1-\alpha(y)]p(x,y)}{\alpha(x)},$$

which can be written as

$$-\sum_{x,y \in A} \alpha(y) \frac{\pi(y)p(y,x) - \pi(x)p(x,y)}{\pi(x)\pi(y)}.$$

Interchanging the roles of $x$ and $y$ above and adding the two expressions, we see that the left-hand side of (1) is

$$\frac{1}{2} \sum_{x,y \in A} [\alpha(x) - \alpha(y)] \frac{\pi(y)p(y,x) - \pi(x)p(x,y)}{\pi(x)\pi(y)},$$

which is 0 by (4).

Finally, we note that (5) is equivalent to (a) and (b) in the statement of the theorem. Statement (b) in the theorem is equivalent to (4). Given this, it suffices to show that statement (a) is equivalent to (2) for all $u$. But the left-hand side of (2) for $u \in C_i$ can be written as

$$(1 - \alpha_i) \sum_j (1 - \alpha_j) \sum_{y \in C_j} [\pi_j p(y,u) - \pi_i p(u,y)].$$

The summands above for $j \neq i$ vanish by (b), so (2) for this $u$ is equivalent to

$$\sum_{y \in C_i} [\pi_i p(y,u) - \pi_i p(u,y)] = 0.$$

But this is statement (a), so the theorem is proved. $\square$

**3. Product stationary measures on $S = Z^d$.** In this section we assume that $S = Z^d$ and $p(x,y) = p(y-x)$ for all $x, y$. We will show later that, under a minimal assumption, the $\pi$ corresponding to a stationary product measure must have an exponential form. [Recall $\pi(x) = \alpha(x)/(1-\alpha(x))$ if the product measure is $\nu_\alpha$.] In the next result, we obtain an easily checked necessary and sufficient condition for stationarity in this case.

THEOREM 2. *Suppose that $\pi(x) = \pi(0)e^{\langle x,v \rangle}$ for all $x$ and some $v \in R^d$. Then, $\nu_\alpha$ is stationary for $\eta_t$ if and only if*

(6) $\qquad p(z) = e^{\langle z,v \rangle} p(-z) \qquad$ *for all $z$ such that $\langle z, v \rangle \neq 0$.*

PROOF. Condition (a) of Theorem 1 says in this case that

$$\sum_{y:\, \langle y,v \rangle = \langle x,v \rangle} [p(x,y) - p(y,x)] = 0.$$

Letting $z = y - x$ and using translation invariance, this becomes

$$\sum_{z:\, \langle z,v \rangle = 0} [p(z) - p(-z)] = 0,$$

which is automatically true. Condition (b) of Theorem 1 says that

$$e^{\langle x,v \rangle} p(y-x) = e^{\langle y,v \rangle} p(x-y) \qquad \text{whenever } e^{\langle x,v \rangle} \neq e^{\langle y,v \rangle}.$$

Again, letting $z = y - x$, we see that, in this case, (6) is equivalent to (b) of Theorem 1. $\square$

REMARK. Theorem 2 makes it easy to construct examples of exclusion processes with a significant number of nontranslation invariant stationary measures of product form, including many that are not reversible. For example, suppose $d = 2$,

$$p(e_1) = p_1, \qquad p(e_2) = p_2, \qquad p(-e_1) = q_1, \qquad p(-e_2) = q_2$$

and $p(z) = 0$ for all other $z$'s. Here $e_1, e_2$ are the standard basis elements in $Z^2$. Assume that $p_1 \neq q_1$ and $p_2 \neq q_2$. Then the stationary product measures with density $1/2$ at the origin are exactly the four corresponding to the following:

$$\pi(x) \equiv 1, \qquad \qquad \pi(x) = (p_1/q_1)^{x^{(1)}},$$
$$\pi(x) = (p_2/q_2)^{x^{(2)}}, \qquad \pi(x) = (p_1/q_1)^{x^{(1)}} (p_2/q_2)^{x^{(2)}},$$

where we have written $x = (x^{(1)}, x^{(2)})$. Only the last of these is reversible.

Next, we will show that, under a weak assumption (that is trivially necessary), $\pi$ must have the exponential form assumed in Theorem 2. Most of the work is contained in the proof of the following statement.

PROPOSITION 1. *Suppose $\alpha$ is such that $\nu_\alpha$ is stationary for $\eta_t$, and $u \in Z^d$ satisfies $p(u) > 0$.*

(a) *If $p(-u) = 0$, then $\pi(x + u) = \pi(x)$ for all $x \in Z^d$.*

(b) *If $p(-u) > 0$, then for each $x \in Z^d$, either*

$$\pi(x + nu) = \pi(x) \quad \text{for all } n \in Z,$$

*or*

$$\pi(x + nu) = \pi(x)\lambda^n \quad \text{for all } n \in Z,$$

*where*

(7) $$\lambda = \frac{p(u)}{p(-u)}.$$

PROOF. Suppose first that $p(-u) = 0$, and fix an $x$. Then the equality in (b) of Theorem 1 cannot hold with $y = x + u$ since one side is zero and the other is not. Therefore, $\pi(x) = \pi(x + u)$. This proves part (a) of the proposition.

The proof of part (b) is significantly harder. Suppose $p(-u) > 0$, and define $\lambda$ as in (7). By Theorem 1(b),

(8) $$\text{for each } x \in Z^d, \pi(x+u)/\pi(x) = 1 \quad \text{or} \quad \lambda.$$

Therefore, for each $x$ and $n$, $\pi(x + nu)/\pi(x)$ is a power of $\lambda$, and the sign of the power is the same as the sign of $n$. Our problem is to show that, for each $x$, the power is identically $n$ or identically 0. To do so, we will need to use part (a) of Theorem 1 as well. Note that the result is immediate if $\lambda = 1$, so we may assume that $\lambda \neq 1$. By interchanging the roles of $u$ and $-u$ if necessary, we may assume that $\lambda > 1$.

As mentioned above, we will need to use the property in Theorem 1(a). This involves sums in which $p(x, y)$ and $p(y, x)$ both appear for $x, y$ with the same value of $\pi$. It will be important to relate these two transition probabilities so that one sum can be related to the other. That is the purpose of the next observations. We will call the two sites of interest $z$ and $w$ instead of $x$ and $y$ to emphasize that, for the moment, they are fixed, rather than being summed over. In (10) and (11) below, we note that, depending on how $z$ and $w$ are placed relative to the values of $\pi$ at these and "neighboring" points, one of the transition probabilities, $p(z, w)$ and $p(w, z)$, is $\lambda$ times the other.

Consider two distinct points $z, w \in Z^d$ for which $\pi(z) = \pi(w)$. If $p(z, w) > 0$ and $p(w, z) = 0$ or vice versa, then by part (a) of this proposition,

(9) $$\pi(z + nu) = \pi(w + nu) \quad \text{for all } n \in Z.$$

Suppose now that $p(z, w) > 0$ and $p(w, z) > 0$. If $\pi(w + u) \neq \pi(w)$, then by (8),

$$\pi(w + u)/\pi(w) = \lambda.$$

If also $\pi(z+u) = \pi(z)$, then $\pi(w+u) \neq \pi(z+u)$, so by Theorem 1(b),
$$p(z,w) = \lambda p(w,z).$$
Similarly, if $\pi(w-u) \neq \pi(w)$ and $\pi(z-u) = \pi(z)$, then
$$p(w,z) = \lambda p(z,w).$$

Now let $c$ be a value of the function $\pi(\cdot)$, and write
$$A = \{z : \pi(z-u) = \lambda^{-1}c, \pi(z) = c, \pi(z+u) = \lambda c\},$$
$$B = \{z : \pi(z-u) = \lambda^{-1}c, \pi(z) = c, \pi(z+u) = c\}$$
and
$$C = \{z : \pi(z-u) = c, \pi(z) = c\}.$$
Note that by (8),
$$A \cup B \cup C = \{z : \pi(z) = c\}.$$
If $w \in A \cup B$ and $z \in C$, (9) is false (for $n = -1$), so $p(z,w)$ and $p(w,z)$ are either both positive or both zero. If they are both positive, then we know from the previous paragraph that $p(w,z) = \lambda p(z,w)$ and, of course, this is also true if they are both zero. So, it is always the case that

(10) $\qquad w \in A \cup B, z \in C \qquad$ implies $p(w,z) = \lambda p(z,w)$.

Similarly,

(11) $\qquad w \in B, z \in A \qquad$ implies $p(w,z) = \lambda p(z,w)$.

By Theorem 1(a), since $\{z : \pi(z) = c\} = A \cup B \cup C$,

(12) $\qquad \displaystyle\sum_{z \in A \cup C}[p(z,w) - p(w,z)] + \sum_{z \in B}[p(z,w) - p(w,z)] = 0, \qquad w \in B.$

Combining (10), (11) and (12) gives

(13) $\qquad (\lambda - 1)\displaystyle\sum_{z \in A \cup C} p(w,z) = \lambda \sum_{z \in B}[p(z,w) - p(w,z)], \qquad w \in B.$

Note that $w \in B$ implies $z = w + u \in C$, so that (13) implies

(14) $\qquad \displaystyle\sum_{z \in B}[p(z,w) - p(w,z)] \geq (1 - \lambda^{-1})p(u) > 0, \qquad w \in B.$

If $d = 1$, we can now complete the proof easily, since in this case, $B$ is finite by (8). The sum on $w \in B$ of the left-hand side of (14) is zero, so the left-hand side itself is zero as well. It follows that $B = \varnothing$. Similarly,
$$\{z : \pi(z-u) = c, \pi(z) = c, \pi(z+u) = \lambda c\} = \varnothing.$$

Together with (8), this completes the proof of the second part of the proposition in this case.

An extra argument is needed in higher dimensions. The objective is still to show that $B = \varnothing$ as a consequence of (14). Once this is done, the proof of the proposition is completed as in the previous paragraph.

To do so, assume $B \neq \varnothing$ and let $D \subset B$ be finite and nonempty. Sum (14) for $w \in D$. The part of the resulting sum on the left that corresponds to $w, z \in D$ vanishes. Therefore, making the change of variable $z = w + x$ and dividing by the cardinality of $D$, one obtains

$$(1 - \lambda^{-1})p(u) \leq \sum_x \left[|p(-x) - p(x)| \frac{\#\{w \in D : w + x \in B \setminus D\}}{|D|}\right].$$

Now let $F \subset Z^d$ be finite. It follows from the above inequality that

$$(1 - \lambda^{-1})p(u) \leq \sum_{x \notin F} |p(-x) - p(x)| + \frac{\#\{(x, w) \in F \times D : w + x \in B \setminus D\}}{|D|}.$$

To get a contradiction, it will suffice to find, for any finite $F$, an increasing sequence $D_n$ of finite, nonempty subsets of $B$ so that

(15) $$\liminf_n \frac{\#\{(x, w) \in F \times D_n : w + x \in B \setminus D_n\}}{|D_n|} = 0.$$

Once this is done, it will follow that the sums on the right-hand sides of the previous two displays can be made arbitrarily small by choosing $F$ sufficiently large.

Let $D_0 = \varnothing$, $D_1$ be any finite nonempty subset of $B$, and define recursively

$$D_{n+1} = \{w + x : x \in F, w \in D_n, w + x \in B\} \cup D_n.$$

If (15) fails for this sequence, then there is an $\varepsilon > 0$ so that

(16) $$\#\{(x, w) \in F \times D_n : w + x \in B \setminus D_n\} \geq \varepsilon |D_n|, \qquad n \geq 1.$$

To show that this is not possible, consider the mapping

$$(x, w) \to x + w.$$

This maps $\{(x, w) \in F \times D_n : w + x \in B \setminus D_n\}$ into $D_{n+1} \setminus D_n$, and is at most $|F|$ to 1. Therefore, by (16),

(17) $$\varepsilon |D_n| \leq |F| |D_{n+1} \setminus D_n|,$$

so that $|D_n|$ grows exponentially rapidly. On the other hand, $D_n$ is contained in the union of $|D_1|$ balls of radius $Kn$, where $K = \max\{\|x\|, x \in F\}$, so $|D_n|$ can grow at most like $K'n^d$ for some constant $K'$. This contradicts the exponential growth, so (15) holds. □

Next we extend the statement of Proposition 1 to several values of $u$.

PROPOSITION 2. *Suppose $\alpha$ is such that $\nu_\alpha$ is stationary for $\eta_t$, and that $p(u_i) > 0$ for $1 \leq i \leq k$. If $p(-u_i) > 0$, set*

$$\lambda_i = \frac{p(u_i)}{p(-u_i)}, \tag{18}$$

*while if $p(-u_i) = 0$, set $\lambda_i = 1$. Then*

$$\pi(n_1 u_1 + \cdots + n_k u_k) = \pi(0) \rho_1^{n_1} \cdots \rho_k^{n_k}, \tag{19}$$

*where for each $i$, $\rho_i = \lambda_i$ or $1$.*

PROOF. We prove this for $k = 2$, since the proof for general $k$ follows by iterating the argument. So, suppose $p(u_1) > 0$ and $p(u_2) > 0$. By Proposition 1, for each $x$,

$$\pi(x + mu_1 + nu_2) = \begin{cases} \pi(x + mu_1)(\rho_m)^n, \\ \pi(x + nu_2)(\sigma_n)^m, \end{cases}$$

where $\rho_m = \lambda_2$ or $1$ for each $m$, and $\sigma_n = \lambda_1$ or $1$ for each $n$. Putting $n = 0$ and then $m = 0$ in the above expression, and using the resulting expressions, we see that

$$\pi(x + mu_1 + nu_2) = \pi(x)(\sigma_0)^m(\rho_m)^n = \pi(x)(\rho_0)^n(\sigma_n)^m$$

for all $m, n$. It follows that

$$(\rho_m/\rho_0)^{1/m} = (\sigma_n/\sigma_0)^{1/n}.$$

Since the left-hand side does not depend on $n$ and the right-hand side does not depend on $m$, this expression is a constant $c$ independent of $m$ and $n$:

$$\rho_m = \rho_0 c^m \quad \text{and} \quad \sigma_n = \sigma_0 c^n.$$

Since $\rho_m$ and $\sigma_n$ only have two possible values, we conclude that $c = 1$, and therefore that $\rho_m$ is independent of $m$ and $\sigma_n$ is independent of $n$. Therefore, for each $x$,

$$\pi(x + mu_1 + nu_2) = \pi(x)\sigma^m \rho^n,$$

where $\rho = \lambda_2$ or $1$, and $\sigma = \lambda_1$ or $1$. (The choice of $\rho$ and $\sigma$ can depend on $x$, but not on $m$ or $n$.) □

We now use Proposition 2 to prove the main result of this section.

THEOREM 3. *Assume that the translation invariant exclusion process with transition probabilities $p(x, y) = p(y - x)$ has the property that there is no proper subgroup of $Z^d$ that contains $P = \{u \in Z^d : p(u) > 0\}$. If $\alpha$ is such that $\nu_\alpha$ is stationary, then there is a $v \in R^d$ so that $\pi(x) = \pi(0) e^{\langle x, v \rangle}$ for all $x$.*

PROOF. Let $u_1, \ldots, u_k$ be elements of $P$ that span $Z^d$. To prove the theorem, it suffices to show that there is a $v \in R^d$ so that

(20) $$\langle u_i, v \rangle = \log \rho_i, \qquad 1 \leq i \leq k,$$

where the quantities $\rho_i$ are the ones that appear in the statement of Proposition 2. By relabelling, we may assume that $u_1, \ldots, u_d$ is a basis for the vector space $Q^d$. Let $v \in R^d$ be the unique solution of the equations (20) for $1 \leq i \leq d$. For $d < i \leq k$, $u_i$ can be written as a linear combination of $u_1, \ldots, u_d$,

$$nu_i = \sum_{j=1}^{d} n_j u_j,$$

where $n$ and $n_1, \ldots, n_d$ are integers. Applying (19) twice gives

$$\pi(nu_i) = \pi(0)\rho_i^n \quad \text{and} \quad \pi(n_1 u_1 + \cdots + n_d u_d) = \pi(0)\rho_1^{n_1} \cdots \rho_d^{n_d}.$$

Therefore,

$$\rho_i^n = \rho_1^{n_1} \cdots \rho_d^{n_d},$$

so that

$$n \log \rho_i = \sum_{j=1}^{d} n_j \log \rho_j = \sum_{j=1}^{d} n_j \langle u_j, v \rangle = n \langle u_i, v \rangle,$$

and (20) holds for this $i$ as well. □

**4. Product stationary measures on trees.** For most of this section, we take $S$ to be the homogeneous tree $T$ in which every vertex has degree 3. (Similar examples with a general degree could be analyzed in a similar way, but we consider this case for concreteness.) Suppose $q, r, s$ are positive, distinct and sum to 1. We will use Theorem 1 to determine all stationary product measures for various exclusion processes on $T$ that are homogeneous and nearest neighbor. We will initially consider $p(x, y)$ so that the corresponding Markov chain has the property that, from each vertex, there are probabilities $q, r, s$ of going out to the three neighbors, and also probabilities $q, r, s$ coming in from the three neighbors. There are three inequivalent homogeneous ways of making these assignments that are given below as Cases 1, 2 and 3. The existence of nonreversible, inhomogeneous (i.e., with nonconstant density) stationary product measures for the process depends on the assignments.

CASE 1. The edge joining $x$ and $y$ is said to be of type

1     if $\{p(x,y), p(y,x)\} = \{q, r\}$,

2     if $\{p(x,y), p(y,x)\} = \{r, s\}$,

3     if $\{p(x,y), p(y,x)\} = \{s, q\}$.

Assign labels $1, 2, 3$ to the edges of $T$ in such a way that for every vertex $x$ the three edges incident to $x$ consist of one of each of the three types. Up to isomorphism, there is only one way to do this. Once the labels are assigned, assign transition probabilities $q, r, s$ in each direction for each edge in such a way that these assignments are consistent with the labels, and for each vertex, there are probabilities $q, r, s$ of moving to the three neighbors. Again, up to isomorphism, there is only one way to do this. (In each case, the meaning of the word isomorphism should be clear. In the first occurrence, e.g., it means that, for any two assignments of labels to edges, there is a 1–1 map from $T$ with one assignment onto $T$ with another assignment that respects the labels.)

For a prescribed value of $\pi(x_0)$ for a given vertex $x_0$, there is only one choice of $\pi$ for which $\nu_\alpha$ is reversible for the exclusion process with this kernel $p(\cdot)$; it is determined by

$$\pi(x)p(x,y) = \pi(y)p(y,x)$$

for all neighbors $x, y$. We will show that the only stationary product measures for this process are the homogeneous measures and the reversible measures.

Take any vertex $x$, and let $x_1, x_2, x_3$ be the three neighbors of $x$, with the indexes chosen so that the edge joining $x$ and $x_i$ is of type $i$. If $\nu_\alpha$ is stationary, then by Theorem 1(a),

$$(q-r)\mathbb{1}_{\{\pi(x_1)=\pi(x)\}} + (r-s)\mathbb{1}_{\{\pi(x_2)=\pi(x)\}} + (s-q)\mathbb{1}_{\{\pi(x_3)=\pi(x)\}} = 0.$$

It follows that, for each $x$, either all or none of the neighbors $x_i$ of $x$ have the property that $\pi(x_i) = \pi(x)$. However, it is clear that which of these situations occurs is independent of $x$, since any two neighbors must be in the same situation. Therefore, either $\pi$ is constant on $T$, or $\pi(x) \neq \pi(y)$ for all neighbors $x, y$. In either case, $\pi$ is determined by its value at one site [by Theorem 1(b) in the latter case]. So, up to scaling, there are two stationary product measures, the first homogeneous and the second reversible. Note that these correspond to the two cases in Theorem 2.1 on page 380 in [13].

CASE 2. The edge joining $x$ and $y$ is said to be of type

1      if $\{p(x,y), p(y,x)\} = \{q, r\}$,

2      if $\{p(x,y), p(y,x)\} = \{s, s\}$.

Assign labels $1, 2$ to the edges of $T$ in such a way that, for every vertex $x$, two of the three edges incident to $x$ are of type 1 and the other is of type 2. Again, up to isomorphism, there is only one way to make the assignment. Assign transition probabilities to the edges in a manner consistent with the assignment of edge types. Each site $x$ has exactly one neighbor $y$ with $p(x,y) = q$, and one neighbor $y$ with $p(x,y) = r$.

To describe all functions $\pi$ corresponding to stationary product measures, call a subset $L$ of vertices of $T$ a line if it is isomorphic to $Z$, and all edges joining vertices in $L$ are of type 1. Then $T$ can be partitioned into lines in a natural manner. Each edge of type 2 joins two lines. Since the transition probabilities are symmetric across edges of type 2, $\pi$ must take the same value at two vertices that are joined by an edge of type 2 [by Theorem 1(b)]. On each line, $\pi$ is either constant or of the form

$$\pi(x_n) = \pi(x_0)(q/r)^n; \tag{21}$$

one sees this by applying Theorem 1(a), together with the isomorphism in the form $L = \{\ldots, x_{-1}, x_0, x_1, \ldots\}$ with the appropriate orientation. So, the most general choice of $\pi$ so that $\nu_\alpha$ is stationary is obtained by (a) fixing the value of $\pi(x)$ at one vertex $x$, and (b) deciding for each line $L$ whether $\pi$ is to be constant on $L$, or of the form (21). The stationary product measure is homogeneous if $\pi$ is constant on all lines, and is reversible if $\pi$ is of the form (21) on all lines. In all other cases, $\nu_\alpha$ is neither homogeneous nor reversible.

CASE 3. The edge joining $x$ and $y$ is said to be of type

$$\begin{array}{ll} 1 & \text{if } \{p(x,y), p(y,x)\} = \{q,q\}, \\ 2 & \text{if } \{p(x,y), p(y,x)\} = \{r,r\}, \\ 3 & \text{if } \{p(x,y), p(y,x)\} = \{s,s\}. \end{array}$$

Make consistent assignments of edge types and transition probabilities as before. In this case, the process is symmetric, so by Theorem 1(b), the only stationary product measures are the homogeneous ones, and they are, of course, reversible.

In all three of the above cases, the homogeneous product measures are stationary, but this is not necessarily the case for other homogeneous transition probabilities on the above tree $T$. For instance, one can choose the probabilities from a vertex to its three neighbors to be $q, r, s$ as before, but with $q = r \neq s$ (and, hence, $s = 1 - 2q$). Think of the tree as branching up, so that each vertex has two edges that go up and one that goes down. The probability $q$ is assigned to each of the two upward edges and the probability $s$ to the downward edge from each vertex. The corresponding Markov chain is not doubly stochastic. So by Theorem 1(a), there are no (nontrivial) homogeneous stationary product measures, although the transition probabilities are themselves homogeneous. This, of course, does not occur for $S = Z^d$, where the Markov chain is doubly stochastic for such transition probabilities.

The exponential growth of the number of neighbors within distance $D$ of a given site in a tree can give rise to a phase transition for the corresponding stationary product measures. Consider, for example, the exclusion process in the preceding paragraph, but restricted to the (inhomogeneous) tree consisting of all vertices of $T$ that are equal to or above a fixed vertex $x_0$ (with particles prevented from leaving this subtree). By Theorem 1, the stationary product measures for the process are reversible and satisfy

$$\pi(x) = (q/(1-2q))^n \pi(x_0) \qquad \text{for } D(x_0, x) = n.$$

When are such stationary product measures extremal stationary measures? Theorem 2.1 in [7] gives, as a necessary and sufficient criterion [in the more general setting of countable $S$ with irreducible $p(\cdot)$], that

$$\sum_{x \in S} \pi(x)/(1+\pi(x))^2 = \infty.$$

Here, this occurs for $q/(1-2q) \in [1/2, 2]$, that is, when $q \in [1/4, 2/5]$. This contrasts with $S = Z$, where none of the inhomogeneous stationary product measures is extremal, and with $S = Z^d$, $d \geq 2$, where all stationary product measures are extremal.

**5. Stationary profile measures—necessary conditions for existence.** In the next two sections we consider the exclusion process $\eta_t$ on $Z^d$ with transition probabilities $p(x, y) = p(y - x)$ that satisfy

$$\sum_x |x| p(x) < \infty.$$

The analogues of blocking or profile measures in higher dimensions for such $p(\cdot)$ should be measures that are invariant under translations in certain directions, and have some specified limiting behavior in directions orthogonal to those directions. In this section we rule out the existence of stationary measures with these properties in certain directions. These results are analogues of the statement in one dimension that, if the particle drift is positive, and a stationary measure has density tending to zero in one direction and to one in the opposite direction, then it must be that the limit is one in the $+\infty$ direction and zero in the $-\infty$ direction.

We begin with some definitions. Fix a nonzero vector $v \in R^d$. If $v \in Z^d$, we want to say that a measure $\mu$ is $v$-*homogeneous* if it is invariant under translations in all directions orthogonal to $v$. Unfortunately, this statement makes no sense for general $v \in R^d$, since translations in directions orthogonal to $v$ will, in general, lead outside of $Z^d$. In most cases, we will consider $v \in Z^d$, but to prepare for situations where $v$ is more general (e.g., Theorem 5 below), we make a definition that is meaningful for all $v \in R^d$, and agrees with the previous statement if $v \in Z^d$. In general, a measure $\mu$ on $\{0, 1\}^{Z^d}$ will be called *v-homogeneous* if it has the following property: For each finite $A \subset Z^d$, there is a continuous function $f_A$ on $R$ so that

$$(22) \qquad \mu\{\eta : \eta \equiv 1 \text{ on } x + A\} = f_A(\langle x, v \rangle), \qquad x \in Z^d.$$

It is easy to check that if $v \in Z^d$, this is equivalent to $\mu$ being invariant under shifts in all directions orthogonal to $v$, as we wanted. Note that the product measures $\nu_\alpha$ considered in Theorem 2 are $v$-homogeneous for general $v \in R^d$, since in that case,

$$\nu_\alpha\{\eta : \eta(x) = 1\} = \frac{\pi(0)}{\pi(0) + e^{-\langle x, v \rangle}}.$$

This supports the view that our definition is a good one in this context. Of course, a measure is $v$-homogeneous if and only if it is $cv$-homogeneous for any $c \in R \setminus \{0\}$, so there is no difference in statements of hypotheses between assuming $v \in Z^d$ and assuming $v \in Q^d$.

The continuity assumption on $f_A$ in the above definition is not an additional requirement if $v \in Z^d$, since then $\langle x, v \rangle$ only takes integer values. On the other hand, for general $v \in R^d$, if we did not require $f_A$ to be continuous, property (22) would contain little, if any, information about $\mu$. The problem is that the map $x \to \langle x, v \rangle$ could easily be one-to-one, and in this case, any function of $x$ is a function of $\langle x, v \rangle$. In remarks following the proofs of Theorems 4 and 6, we will indicate where the continuity property would arise if we were to take $v \in R^d$ instead of $v \in Z^d$.

Note that if $v \in Z^d$, for example, (22) does not determine the values of $f_A(s)$ for all $s \in R$. However, by choosing natural versions of $f_A$, we may assume for general $v \in R^d$ that

$$(23) \qquad f_{x+A}(s) = f_A(s + \langle x, v \rangle) \qquad \text{for } x \in Z^d \text{ and } s \in R.$$

An important property of $v$-homogeneity is that the class of $v$-homogeneous measures is preserved by the evolution of the exclusion process. This is clear if $v \in Z^d$. To see it in general, let $S(t)$ be the semigroup for the exclusion process, and for finite $A \subset Z^d$, let $\chi_A(\eta)$ be the indicator function of $\{\eta : \eta \equiv 1 \text{ on } A\}$. Then $S(t)\chi_A(\eta)$ is continuous in $\eta$, so it can be uniformly approximated by functions of the form $\sum_{i=1}^n c_i \chi_{B_i}$. Furthermore, since the exclusion process is translation invariant, $S(t)\chi_{x+A}$ is then automatically uniformly approximated by $\sum_{i=1}^n c_i \chi_{x+B_i}$, where the uniformity applies to both $\eta$ and $x$. Therefore, if $\mu$ is $v$-homogeneous,

$$\mu S(t)\{\eta : \eta \equiv 1 \text{ on } x + A\} = \int S(t) \chi_{x+A} \, d\mu$$

is uniformly (in $x$) approximated by

$$\sum_{i=1}^n c_i \int \chi_{x+B_i} \, d\mu = \sum_{i=1}^n c_i f_{B_i}(\langle x, v \rangle).$$

It follows that $\mu S(t)\{\eta : \eta \equiv 1 \text{ on } x + A\}$ is a continuous function of $\langle x, v \rangle$ as required.

Unfortunately, the class of $v$-homogeneous measures is not preserved under weak convergence for general $v$ because of the continuity assumption on $f_A$. To see this, note that a product measure $\nu_\alpha$ is $v$-homogeneous if and only if $\alpha(x)$ is a continuous function of $\langle x, v \rangle$. Suppose then that $f_n$ is a sequence of continuous functions on $R$ that converges pointwise to a function $f$ that has a jump discontinuity. Let $\alpha_n(x) = f_n(\langle x, v \rangle)$ and $\alpha(x) = f(\langle x, v \rangle)$. If the coordinates of $v$ are linearly independent over the rationals (so that $\{\langle x, v \rangle : x \in Z^d\}$ is dense in $R$), then $\nu_{\alpha_n}$ is $v$-homogeneous, $\nu_\alpha$ is not, yet $\nu_{\alpha_n}$ converges weakly to $\nu_\alpha$. As a consequence of this, we take $v \in Z^d$ in most of our results.

A $v$-homogeneous measure that is asymptotically equal to $\delta_0$ in the $-v$ direction and asymptotically equal to $\delta_1$ in the $v$ direction, in the sense that

$$(24) \qquad \lim_{s \to -\infty} f_A(s) = 0 \quad \text{and} \quad \lim_{s \to \infty} f_A(s) = 1, \qquad A \neq \varnothing,$$

will be called a *v-profile measure*. Our main objective in this section is to demonstrate the nonexistence of $v$-profile stationary measures under suitable assumptions on $p(\cdot)$ and $v$. In the next section we will provide some information about the structure of the set of $v$-profile stationary measures in case of existence.

The results of Section 3 give some examples of existence of $v$-profile stationary measures. For instance, in the remark following Theorem 2, we considered a general nearest neighbor process on $Z^2$. Suppose for concreteness that $p_1 > q_1$ and $p_2 > q_2$. Then, the four stationary product measures $\mu$ described there are, respectively: (i) homogeneous, (ii) $v$-profile with $v = (1, 0)$, (iii) $v$-profile with $v = (0, 1)$, and (iv) $v$-profile with

$$v = \left( \log \frac{p_1}{q_1}, \log \frac{p_2}{q_2} \right).$$

This example illustrates why we do not want to restrict consideration to $v \in Z^d$ (or, equivalently, to $v \in Q^d$). Even if the $p_i$ and $q_i$ are all rational, the above $v$ need not be a multiple of a vector in $Z^2$.

The proof of Theorem 4 below (and that of later results) uses the coupled process $(\eta_t, \zeta_t)$ which has proved so useful in the study of the exclusion process. Each of the two coordinates is assumed to be a copy of the exclusion process, and the coupling is defined by saying that particles in the two processes move together as much as possible. This coupling has the property that sites at which $\eta_t(x) \neq \zeta_t(x)$ can move and disappear, but they cannot be created.

To describe the joint process in more detail, call a site $x$ at which $\eta(x) \neq \zeta(x)$ a *discrepancy*. There are two types of discrepancies, according to whether $\eta(x) = 1$, $\zeta(x) = 0$ or $\eta(x) = 0$, $\zeta(x) = 1$. Say that a site $x$ with no discrepancy is of *type* 1 if $\eta(x) = \zeta(x) = 1$ and of *type* 0 if $\eta(x) = \zeta(x) = 0$. Consider now the possible transitions for the coupled process involving two sites $x \neq y$. If $x$ and $y$ are both of type 1 or both of type 0, nothing happens. If $x$ is of type 1 and $y$ is of type 0, the particles at $x$ in both configurations move together to $y$ at rate $p(y - x)$. If there are discrepancies at both $x$ and $y$, nothing happens if they are of the same type. If they are of opposite type, the particles in the two configurations move independently to the other site at the appropriate rates [one at rate $p(y-x)$ and the other at rate $p(x-y)$]; after the transition, one of the sites $x$ and $y$ is of type 1 and the other is of type 0, that is, two discrepancies of the process have disappeared. Finally, if there is a discrepancy at $x$, and $y$ is of type 0 or 1, then the particle in the configuration for which the values at $x$ and $y$ differ moves to the other site at the appropriate rate. For more on this coupling and its application to the exclusion process, see Section 2 of Chapter VIII in [13].

The concepts $v$-homogeneous and $v$-profile have natural analogues for the measures on $\{0,1\}^{Z^d} \times \{0,1\}^{Z^d}$ that arise as the measure of the coupled process. To define the first, the left-hand side of (22) is replaced by any cylinder probability, while for the second, one can simply say that $\mu$ on $\{0,1\}^{Z^d} \times \{0,1\}^{Z^d}$ is $v$-profile if it is $v$-homogeneous and each of its two marginals is $v$-profile.

Our next result, Theorem 4, asserts that existence of $v$-profile stationary measures is only possible if $v$ and the mean vector of $p(\cdot)$ have nonnegative inner product. Unfortunately, the proof requires that $v$ be in $Z^d$. The problem that arises when $v$ is general is closely connected to the continuity assumption that we made in the definition of $v$-homogeneity. (For more on this point, see the remark following the proof.)

The proof of Theorem 4 is based on the idea of the proof of Proposition 3.4 in [1]. It involves a limiting argument that requires some domination. The following two lemmas provide the domination and main limiting statement that are needed. Neither involves the exclusion process or the kernel $p(\cdot)$. To begin, let $v_1, \ldots, v_d \in Z^d$ be an orthogonal basis for $R^d$, with $v_1 = v$. Such a choice exists since the Gram–Schmidt orthogonalization procedure produces vectors with rational coordinates if the original vectors have rational coordinates. For a positive integer $n$, let

$$B_n = \{x \in Z^d : |\langle x, v_i \rangle| \leq n, 1 \leq i \leq d\}.$$

The statements of the two lemmas should be intuitively clear. (The precise constant in Lemma 1 is not important.) The reader should simply think of replacing the sums over points in $Z^d$ with integrals over $R^d$ in Lemma 1, for example.

LEMMA 1.
$$\sum_u |\mathbb{1}_{B_n}(u) - \mathbb{1}_{B_n}(u-z)| \leq 2(2n+1)^{d-1} \sum_{j=1}^d |\langle z, v_j \rangle|.$$

PROOF. Write the left-hand side above as

(25)
$$\sum_u |\mathbb{1}_{B_n}(u) - \mathbb{1}_{B_n}(u-z)|$$
$$= \#(u \in B_n : u - z \notin B_n) + \#(u \notin B_n : u - z \in B_n).$$

The first term on the right-hand side is at most

(26)
$$\sum_{j=1}^d \#(u : |\langle u, v_i \rangle| \leq n \; \forall 1 \leq i \leq d, |\langle u, v_j \rangle - \langle z, v_j \rangle| > n).$$

The mapping $T$ from $Z^d$ to $Z^d$ defined by $Tu = (\langle u, v_1 \rangle, \ldots, \langle u, v_d \rangle)$ is one-to-one (though not onto). Write $u = (u^{(1)}, \ldots, u^{(d)})$ in Cartesian coordinates, and define
$$C = \{u : |u^{(i)}| \leq n \; \forall 1 \leq i \leq d, \; |u^{(j)} - \langle z, v_j \rangle| > n\}$$
for fixed $j$ and $n$. Then the $j$th summand in (26) is the number of elements in $T^{-1}C$. This is at most the number of elements in $C$, which is, in turn, at most $(2n+1)^{d-1} |\langle z, v_j \rangle|$. Summing on $j$ and applying the same estimate to the second term on the right-hand side of (25), we see that (25) is bounded above by

(27)
$$2(2n+1)^{d-1} \sum_{j=1}^d |\langle z, v_j \rangle|,$$

which completes the proof of the lemma. □

LEMMA 2. *There is a subsequence $n_k$ and a positive constant $K$ so that*

(28)
$$\lim_{k \to \infty} n_k^{1-d} \#(u \in B_{n_k} : \langle u, v \rangle = r) = K$$

*for every $r \in \Gamma = \{\langle u, v \rangle : u \in Z^d\} = aZ$ for some $a > 0$.*

PROOF. We first show that $n^{1-d} \#(u \in B_n : \langle u, v \rangle = 0)$ is bounded away from 0 and $\infty$. As already observed, the mapping $T$ defined in the proof of Lemma 1 is one-to-one, so
$$\#(u \in B_n : \langle u, v \rangle = 0) = \#(u \in Z^d : (Tu)^{(1)} = 0, |(Tu)^{(j)}| \leq n \; \forall 2 \leq j \leq d)$$
$$\leq \#(u \in Z^d : u^{(1)} = 0, |u^{(j)}| \leq n \; \forall 2 \leq j \leq d)$$
$$= (2n+1)^{d-1}.$$

This gives the upper bound. For the lower bound, note that $u = k_2 v_2 + \cdots + k_d v_d \in B_n$ if $k_i \in Z$ and $|k_i| \leq n|v_i|^{-2}$ for $2 \leq i \leq d$, so

$$\#(u \in B_n : \langle u, v \rangle = 0) \geq \prod_{i=2}^{d}(2n|v_i|^{-2} - 1).$$

We may now choose the subsequence $n_k$ so that the limit in (28) exists and is positive when $r = 0$.

To complete the proof of (28), it now suffices to prove that, for $r \in \Gamma$,

$$\lim_{n \to \infty} n^{1-d}[\#(u \in B_n : \langle u, v \rangle = r) - \#(u \in B_n : \langle u, v \rangle = 0)] = 0$$

(which is true along the full sequence). Since $r \in \Gamma$, there is a $w \in Z^d$ so that $\langle w, v \rangle = r$. Therefore,

$$\#(u \in B_n : \langle u, v \rangle = r) = \#(u \in B_n : \langle u - w, v \rangle = 0)$$
$$= \#(u \in B_n - w : \langle u, v \rangle = 0),$$

from which it follows that

$$\#(u \in B_n : \langle u, v \rangle = r) - \#(u \in B_n : \langle u, v \rangle = 0) = O(n^{d-2}),$$

as required, by the argument that led to (27). The only significant difference is that there is an additional constraint in the definition of the analogous $C$—the constraint that $u^{(1)} = 0$. This accounts for the reduction in the power of $n$ from $d-1$ to $d-2$. □

THEOREM 4. *Take $v \in Z^d \setminus \{0\}$ and suppose that $\mu$ is a $v$-profile measure that is stationary for the exclusion process. Then, $\langle \sum_x xp(x), v \rangle \geq 0$. If $p(\cdot)$ is irreducible, the inequality is strict.*

PROOF. Let $\nu$ be a coupling measure that is $v$-homogeneous, is stationary for the coupled process introduced above, and has first and second coordinates with measures $\mu$ and $\nu_{1/2}$, respectively. Such a measure exists. (See page 383 in [13], where this statement is proved without the part about $v$-homogeneity; the same proof applies to the $v$-homogeneous context because the property of $v$-homogeneity is preserved by the evolution of the exclusion process and by its coupled version.) The net rate at which the total number of discrepancies in $B_n$ changes is zero, since the process is in equilibrium. On the other hand, transitions involving discrepancies of opposite type can only lower the number of discrepancies in $B_n$. Since $\mu$ is $v$-profile, there must be discrepancies of opposite type in $B_n$ if $n$ is sufficiently large. In the irreducible case, they will cancel each other out at a rate that is at least of order $n^{d-1}$. Therefore, the net rate at which discrepancies enter $B_n$ is nonnegative in general, and is at least $\varepsilon n^{d-1}$ in the irreducible case. Explicitly,

$$
\begin{aligned}
&\sum_{x \in B_n, y \notin B_n} p(x-y)\nu\{(\eta,\zeta) : \eta(y) \neq \zeta(y), \eta(x) = \zeta(x) = 0\} \\
&+ \sum_{x \in B_n, y \notin B_n} p(y-x)\nu\{(\eta,\zeta) : \eta(y) \neq \zeta(y), \eta(x) = \zeta(x) = 1\} \\
&\geq \sum_{x \in B_n, y \notin B_n} p(x-y)\nu\{(\eta,\zeta) : \eta(x) \neq \zeta(x), \eta(y) = \zeta(y) = 1\} \\
&+ \sum_{x \in B_n, y \notin B_n} p(y-x)\nu\{(\eta,\zeta) : \eta(x) \neq \zeta(x), \\
&\hspace{6cm} \eta(y) = \zeta(y) = 0\} + \varepsilon n^{d-1},
\end{aligned}
$$
(29)

where $\varepsilon \geq 0$ in general, and $\varepsilon > 0$ in the irreducible case.

To take advantage of the $v$-homogeneity assumption, define functions $g_z$ and $h_z$ by

$$
\begin{aligned}
g_z(\langle x, v \rangle) &= \nu\{(\eta,\zeta) : \eta(x+z) \neq \zeta(x+z), \eta(x) = \zeta(x) = 0\}, \\
h_z(\langle x, v \rangle) &= \nu\{(\eta,\zeta) : \eta(x+z) \neq \zeta(x+z), \eta(x) = \zeta(x) = 1\}.
\end{aligned}
$$

Then, (29) becomes

$$
\begin{aligned}
&\sum_{x \in B_n, y \notin B_n} p(x-y)g_{y-x}(\langle x, v \rangle) + \sum_{x \in B_n, y \notin B_n} p(y-x)h_{y-x}(\langle x, v \rangle) \\
&\geq \sum_{x \in B_n, y \notin B_n} p(x-y)h_{x-y}(\langle y, v \rangle) \\
&+ \sum_{x \in B_n, y \notin B_n} p(y-x)g_{x-y}(\langle y, v \rangle) + \varepsilon n^{d-1}.
\end{aligned}
$$

Making the changes of variables $z = x - y$ or $z = y - x$ and $u = x$ or $u = y$ in the four sums, and noting that, for any sets $B$ and $C$,

$$
\mathbb{1}_{\{u \in B, u \notin C\}} - \mathbb{1}_{\{u \in C, u \notin B\}} = \mathbb{1}_B(u) - \mathbb{1}_C(u),
$$

this leads to

$$
\begin{aligned}
&\sum_z p(z) \sum_u g_{-z}(\langle u, v \rangle)[\mathbb{1}_{B_n}(u) - \mathbb{1}_{B_n}(u-z)] \\
&\geq \sum_z p(z) \sum_u h_z(\langle u, v \rangle)[\mathbb{1}_{B_n}(u+z) - \mathbb{1}_{B_n}(u)] + \varepsilon n^{d-1}.
\end{aligned}
$$
(30)

Next, we will divide (30) by $n^{d-1}$ and pass to the limit (along a subsequence) using the dominated convergence theorem. The domination is provided by Lemma 1. Using this and Lemma 2, we will prove that, for any function $g$ on $\Gamma$ that satisfies

$$\lim_{l \to -\infty} g(l) = 1 \quad \text{and} \quad \lim_{l \to \infty} g(l) = 0,$$

(31) $$\lim_{k \to \infty} n_k^{1-d} \sum_u g(\langle u, v \rangle)[\mathbb{1}_{B_{n_k}}(u) - \mathbb{1}_{B_{n_k}}(u - z)] = K\langle z, v \rangle / a.$$

Any such function $g$ can be uniformly approximated by functions of the form

$$g(l) = \sum_{j=1}^m c_j \mathbb{1}_{\{l \leq l_j\}},$$

where each $l_j \in \Gamma$ and $\sum_{j=1}^m c_j = 1$. By Lemma 1, it suffices to prove (31) for a function $g$ of the form $g(l) = \mathbb{1}_{\{l \leq r\}}$, where $r \in \Gamma$. In this case, if $\langle z, v \rangle > 0$, the left-hand side of (31) becomes

(32) $$\lim_{k \to \infty} n_k^{1-d} \#(u \in B_{n_k} : r - \langle z, v \rangle < \langle u, v \rangle \leq r),$$

so that (31) follows from (28), since there are $\langle z, v \rangle / a$ elements of $\Gamma$ in the interval $(r - \langle z, v \rangle, r]$. A similar argument gives (31) when $\langle z, v \rangle < 0$.

Note that by (24),

$$\lim_{l \to -\infty} g_z(l) = \tfrac{1}{4}, \qquad \lim_{l \to \infty} g_z(l) = 0,$$

$$\lim_{l \to -\infty} h_z(l) = 0, \qquad \lim_{l \to \infty} h_z(l) = \tfrac{1}{4}$$

for every $z$. Therefore, we may apply (31) to $g(l) = 4g_z(l)$ for any $z$. The conclusion is that the limit of the left-hand side of (30) divided by $n^{d-1}$, along the sequence $n_k$, is $\frac{K}{4a} \sum_z p(z)\langle z, v \rangle$. Similarly, applying (31) to $g(l) = 1 - 4h_z(l)$, it follows that the limit of the right-hand side of (30) divided by $n^{d-1}$, along the sequence $n_k$, is $\varepsilon - \frac{K}{4a} \sum_z p(z)\langle z, v \rangle$. Therefore, (30) implies that $\sum_z p(z)\langle z, v \rangle \geq 2a\varepsilon/K$, as required. $\square$

REMARK. As mentioned earlier, the main difficulty in extending Theorem 4 to general $v$ is related to the continuity assumption on $f_A$ in the definition of $v$-homogeneity. In the proof, we used a $v$-homogeneous coupling measure with prescribed $v$-homogeneous marginals. If $v$ were general, the continuity assumption would be needed in the treatment of (31) where we approximated $g$ by step functions. The issue is the existence of a coupling measure that satisfies the continuity property. In the proofs of analogous results in [13] (see pages 143 and 383), one starts the coupled process with the measure that is the product of the two given marginals, and passes to a Cesaro limit of the measure as $t \to \infty$. What is to guarantee that this limiting measure satisfies the needed continuity property? Perhaps a stationary measure for the marginal or coupled process automatically satisfies $v$-homogeneity with the continuity property. We do not know whether this is always the case.

Under a partial symmetry assumption, it is not difficult to prove a similar result that allows $v \in R^d$:

THEOREM 5. *Take $v \in R^d$ and suppose $p(z) = p(-z)$ for all $z$ for which $\langle z, v \rangle \neq 0$. If $\mu$ is $v$-homogeneous and stationary for the exclusion process, then $\mu\{\eta : \eta(x) = 1\}$ is constant on the smallest subgroup of $Z^d$ containing $P = \{x \in Z^d : p(x) > 0\}$. Therefore, there are no $v$-profile stationary measures for the process in the irreducible case.*

Examples to which this result can be applied are the nearest neighbor exclusion processes in which the transition probabilities are symmetric in all directions but one, say $e_1$, and with a drift in the direction of $e_1$. The conclusion in these examples is that there cannot be any $e_k$-homogeneous stationary measure $\mu$ with nonconstant density if $k > 1$.

PROOF OF THEOREM 5. Applying the generator of the exclusion process to the function $\eta \to \eta(x)$, integrating with respect to $\mu$, and using stationarity, gives, for each $x \in Z^d$,

$$(33) \quad \sum_y p(y-x)\mu\{\eta : \eta(x) = 1, \eta(y) = 0\} = \sum_y p(x-y)\mu\{\eta : \eta(x) = 0, \eta(y) = 1\}.$$

Making the change of variables $z = y - x$ in the sums and using (22) yields

$$(34) \quad \sum_z p(z)[f(\langle x, v \rangle) - f_z(\langle x, v \rangle)] = \sum_z p(-z)[f(\langle x + z, v \rangle) - f_z(\langle x, v \rangle)],$$

where we have used the shorthand

$$f(s) = f_{\{0\}}(s) \quad \text{and} \quad f_z(s) = f_{\{0,z\}}(s).$$

Rewriting (34) gives

$$(35) \quad \sum_z p(z)[f(s) - f_z(s)] = \sum_z p(z)[f(s - \langle z, v \rangle) - f_z(s - \langle z, v \rangle)]$$

for all $s \in \Gamma = \{\langle x, v \rangle, x \in Z^d\}$. In obtaining (35), we have used the relation

$$(36) \quad f_{-z}(s) = f_{\{0,-z\}}(s) = f_{-z+\{0,z\}}(s) = f_{\{0,z\}}(s - \langle z, v \rangle) = f_z(s - \langle z, v \rangle),$$

where we have used (23) in the third equality. So far, we have used practically none of the assumptions of the theorem. We point this out because (35) and (36) will be used later in the proof of Theorem 7.

Now, by the symmetry assumption on $p(\cdot)$, the terms in (35) involving $f_z$ cancel out, since

$$
\begin{aligned}
\sum_z p(z)[f_z(s - \langle z, v \rangle) - f_z(s)] &= \sum_{z:\langle z,v \rangle \neq 0} p(z)[f_{-z}(s) - f_z(s)] \\
&= \sum_{z:\langle z,v \rangle \neq 0} p(-z)[f_{-z}(s) - f_z(s)] \\
&= \sum_{z:\langle z,v \rangle \neq 0} p(z)[f_z(s) - f_{-z}(s)].
\end{aligned}
\tag{37}
$$

We have used (36) in the first step, the symmetry assumption in the second and replacement of $z$ by $-z$ in the third. But the second and fourth expressions in (37) are negatives of one another, so they must be zero.

Using the fact that the sums in (37) vanish, (35) leads to

$$f(s) = \sum_z p(z) f(s - \langle z, v \rangle), \qquad s \in \Gamma.$$

This says that $g(x) = f(\langle x, v \rangle)$ is a bounded harmonic function for the random walk on $Z^d$ with probabilities $p(-z)$ of moving from $x$ to $x + z$. It follows that $f$ is a constant on the irreducible classes for this random walk by the Choquet–Deny theorem—see [3], for example. $\square$

**6. Properties of stationary profile measures.** In this section we study exclusion processes on $Z^d$, $d \geq 2$, whose transition probabilities are translation invariant and have a finite mean. Assume, in addition, that the process is irreducible. We begin by using techniques similar to those in Sections 2 and 3 of Chapter VIII in [13] to determine the structure of $v$-profile stationary measures when they exist. Throughout this section, $v$ will be a fixed nonzero element of $Z^d$.

Theorem 6 below has two statements. We will prove the first in detail, but will only sketch the proof of the second part. Three applications of Theorem 6 follow its proof. In the statement below, the inequality $\mu_1 \leq \mu_2$ refers to stochastic monotonicity. See page 71 in [13] for the definition. Recall that an equivalent statement is that $\mu_1$ and $\mu_2$ can be coupled so that any site at which the $\mu_1$ configuration has a 1 must also have a 1 in the $\mu_2$ configuration.

THEOREM 6. *Suppose that $\mu_1$ and $\mu_2$ are extremal stationary measures for the exclusion process.*

(a) *If $\mu_1$ and $\mu_2$ are $v$-profile, then either $\mu_1 \leq \mu_2$ or $\mu_2 \leq \mu_1$.*
(b) *If $\mu_1$ and $\mu_2$ are $v$-homogeneous and $\langle \sum_z z p(z), v \rangle = 0$, then either $\mu_1 \leq \mu_2$ or $\mu_2 \leq \mu_1$.*

PROOF. (a) It is enough to prove that every extremal stationary $v$-profile measure $\nu$ for the coupled process concentrates on configurations $(\eta, \zeta)$ that have no discrepancies of opposite type, that is, that each such $\nu$ satisfies

$$\nu\{(\eta, \zeta) : \eta(x) = \zeta(y) \neq \eta(y) = \zeta(x)\} = 0$$

for all $x, y \in Z^d$. To do so, define $B_{m,n}$ for $1 \leq m \leq n$ slightly more generally than the corresponding sets $B_n$ used in the proof of Theorem 4:

$$B_{m,n} = \{x \in Z^d : |\langle x, v_1 \rangle| \leq m, |\langle x, v_i \rangle| \leq n \text{ for } 2 \leq i \leq d\},$$

where $v_1, \ldots, v_d \in Z^d$ are orthogonal and $v_1 = v$. [Later in the proof, we will take $m, n \to \infty$ with $m = o(n)$.] We will consider the contributions to the expected (with respect to $\nu$) net rate at which discrepancies in $B_{m,n}$ disappear. The expected total net rate is zero because the coupled process is in equilibrium.

If discrepancies of opposite types appear with positive $\nu$ probability, then they will disappear in $B_{m,n}$ at a rate that is at least of order $n^{d-1}$, since (by $v$-homogeneity) there will be at least this many pairs of sites $x, y$ in $B_{m,n}$ with the same difference $x - y$ and the same (by irreducibility) nonzero probability of having a discrepancy of one type at $x$ and of the other type at $y$. Therefore, it suffices to show that the net rate at which discrepancies enter $B_{m,n}$ is $o(n^{d-1})$.

The idea is that discrepancies can enter $B_{m,n}$ across two boundary faces on which the density of particles is close to 0 in one case and close to 1 in the other (by the $v$-profile assumption), or across the other boundary faces. In the first case, there are roughly $n^{d-1}$ locations to consider, and each makes a contribution to the rate of entry of discrepancies that is $o(1)$, so the total contribution is $o(n^{d-1})$. In the second case, the number of locations is of order $mn^{d-2}$, which is $o(n^{d-1})$ if we choose $m = o(n)$. These correspond to cases (B) and (A) below, respectively.

To carry out the details, consider sites $x \in B_{m,n}, y \notin B_{m,n}$, and the expected rate at which a discrepancy moves from $x$ to $y$ or from $y$ to $x$. For example, the expected rate at which a discrepancy moves from $x$ to $y$ is

$$p(y - x)\nu\{(\eta, \zeta) : \eta(x) \neq \zeta(x), \eta(y) = \zeta(y) = 0\}$$
$$+ p(x - y)\nu\{(\eta, \zeta) : \eta(x) \neq \zeta(x), \eta(y) = \zeta(y) = 1\},$$

with a similar expression for the expected rate at which a discrepancy moves from $y$ to $x$. Sums of these expected rates over appropriate choices of $x$ and $y$ will be bounded in the various cases.

Since $y \notin B_{m,n}$, either $|\langle y, v_1 \rangle| > m$ or $|\langle y, v_i \rangle| > n$ for some $2 \leq i \leq d$. There are two cases to consider:

(A) $|\langle y, v_i \rangle| > n$ for some $i > 1$.

(B) $|\langle y, v_1 \rangle| > m$ and $|\langle y, v_i \rangle| \leq n$ for all $i \neq 1$.

The contributions to the expected rate of motion of discrepancies involving all $x, y$ in case (A) (corresponding to motion from points $x$ inside the box to $y$ outside the box or vice versa) for a given $i$ are bounded above by (letting $z = y - x$)

$$\sum_z [p(z) + p(-z)] \#(x : |\langle x, v_1 \rangle| \leq m, |\langle x, v_j \rangle| \leq n \ \forall j > 1, |\langle x + z, v_i \rangle| > n).$$

This sum is at most

$$Cmn^{d-2} \sum_z |z| p(z)$$

for some constant $C$. This is $o(n^{d-1})$, provided we let $m, n \to \infty$ with $m = o(n)$.

Next, we consider the total contribution to the expected rate at which discrepancies enter $B_{m,n}$ for $x$ and $y$ in case (B). Letting $z = y - x$, we see that this total expected rate is bounded above by

$$\sum_z [p(z) + p(-z)] \sum_{\substack{x \in B_{m,n} \\ |\langle x, v_1 \rangle + \langle z, v_1 \rangle| > m}} \nu\{(\eta, \zeta) : \eta(x + z) \neq \zeta(x + z)\}.$$

The number of terms in the inner sum is bounded by a constant multiple of $|z| n^{d-1}$, uniformly in $n, m$ and $z$. Therefore, by the dominated convergence theorem and the first moment assumption on $p(\cdot)$, the above expression will be $o(n^{d-1})$, provided that, for each $z$, the ratio of the inner sum to the number of summands in the inner sum tends to zero as $m, n \to \infty$. Consider the summands for which $\langle x, v_1 \rangle + \langle z, v_1 \rangle > m$, for example—those for which $\langle x, v_1 \rangle + \langle z, v_1 \rangle < -m$ are handled in a similar manner. For these summands, we have $m - \langle z, v_1 \rangle < \langle x, v_1 \rangle \leq m$, and

$$\nu\{(\eta, \zeta) : \eta(x + z) \neq \zeta(x + z)\} \leq \mu_1\{\eta : \eta(x + z) = 0\} + \mu_2\{\eta : \eta(x + z) = 0\}$$
$$= 2 - f_{\{z\},1}(\langle x, v \rangle) - f_{\{z\},2}(\langle x, v \rangle),$$

where we have used $f_{A,i}$ for $i = 1$ and 2 to denote the function defined in (22) corresponding to the measure $\mu_i$. Since $\mu_1$ and $\mu_2$ are $v$-profile measures, the right-hand side above tends to zero as $m \to \infty$. (Recall that $v = v_1$.) This completes the consideration of case (B).

*Sketch of proof of* (b). The only part of the above argument that uses the assumption that the measures $\mu_i$ are $v$-profile, as opposed to $v$-homogeneous, is the treatment of case (B). In that context, we used the fact that the density of discrepancies near the faces of $B_{m,n}$ that are orthogonal to $v$ is small when $m$ is large because the density of 0's (for the face in the $+v$ direction) or 1's (for the face in the $-v$ direction) is small. The point is that, at a discrepancy, one coordinate must be 0 and the other must be 1. If we only assume that the measures are $v$-homogeneous, the net expected rate at which discrepancies cross the faces of $B_{m,n}$ orthogonal to $v$ has to be shown to be small for a different reason. This reason is, of course, related to the assumption that the mean vector for $p(\cdot)$ is orthogonal to $v$. The argument parallels the proof of Theorem 3.14 on page 391 in [13].

Here is the main idea. In a Cesaro sense (with respect to shifts in the $\pm v$ directions), the coupling measure $\nu$, near the faces of $B_{m,n}$ that are orthogonal to $v$, is nearly translation invariant. It is translation invariant in the directions orthogonal to $v$ by the assumption of $v$-homogeneity, and in the $v$ direction as a result of the Cesaro averaging. So, in the limit as $m \to \infty$, the measure $\nu$ is both stationary for the coupled process and translation invariant. This implies that in this limit, the marginals $\mu_1$ and $\mu_2$ are stationary for the exclusion process and translation invariant. But this implies that they are exchangeable by Theorem 3.9 on page 388 in [13]. Furthermore, they are coupled together so that there are no discrepancies of opposite type.

By decomposing the limiting coupling measure according to the type of discrepancy that occurs, we find that we are essentially in the following situation. We can assume that only one type of discrepancy occurs, so that the measures $\mu_1$ and $\mu_2$ are exchangeable, and are coupled by a measure $\nu$ that satisfies

$$\nu\{(\eta,\zeta) : \eta(x) = 1, \zeta(x) = 0\} = 0$$

for all $x$. We will use a suggestive notation for probabilities related to $\nu$. For example, we will write $\nu\{(\eta,\zeta) : \eta(x) = 0, \eta(y) = \zeta(x) = \zeta(y) = 1\}$ as

$$\nu \begin{pmatrix} 1 & 1 \\ 0 & 1 \end{pmatrix}.$$

Entries on the left correspond to site $x$, while entries on the bottom correspond to configuration $\eta$. The total expected net rate at which discrepancies of this type enter $B_{m,n}$ across the face orthogonal to $v$ in the $+v$ direction in case (B) is then

$$\sum_{\substack{x \in B_{m,n}, \langle y, v_1 \rangle > m, \\ |\langle y, v_j \rangle| \leq n \; \forall j \neq 1}} \left\{ p(x-y) \left[ \nu \begin{pmatrix} 0 & 1 \\ 0 & 0 \end{pmatrix} - \nu \begin{pmatrix} 1 & 1 \\ 0 & 1 \end{pmatrix} \right] \right.$$

$$\left. + p(y-x) \left[ \nu \begin{pmatrix} 1 & 1 \\ 1 & 0 \end{pmatrix} - \nu \begin{pmatrix} 1 & 0 \\ 0 & 0 \end{pmatrix} \right] \right\}.$$

Using the fact that $\eta \leq \zeta$ a.s. with respect to $\nu$, the differences that appear above can be written as

$$\nu \begin{pmatrix} 0 & 1 \\ 0 & 0 \end{pmatrix} - \nu \begin{pmatrix} 1 & 1 \\ 0 & 1 \end{pmatrix} = \mu_2\{\zeta : \zeta(x) = 0, \zeta(y) = 1\} - \mu_1\{\eta : \eta(x) = 0, \eta(y) = 1\}$$

and
$$\nu\begin{pmatrix}1 & 1\\ 1 & 0\end{pmatrix} - \nu\begin{pmatrix}1 & 0\\ 0 & 0\end{pmatrix} = \mu_1\{\eta:\eta(x)=1,\eta(y)=0\} - \mu_2\{\zeta:\zeta(x)=1,\zeta(y)=0\}.$$

By exchangeability of the marginals, the right-hand sides above are independent of $x$ and $y$ and are negatives of one another, so we will call them $A$ and $-A$, respectively. It follows that the total expected net rate at which discrepancies enter $B_{m,n}$, across the face orthogonal to $v$ in the $+v$ direction is, in case (B),

$$A \sum_{\substack{x \in B_{m,n}, \langle y, v_1\rangle > m,\\ |\langle y, v_j\rangle| \le n \ \forall j \ne 1}} [p(x-y) - p(y-x)].$$

Except for the factor of $A$, this can be written as

$$\sum_z p(z)[\#(x \in B_{m,n} : \langle x, v_1\rangle > m + \langle z, v_1\rangle, |\langle x, v_j\rangle - \langle z, v_j\rangle| \le n \ \forall j \ne 1)$$
$$- \#(x \in B_{m,n} : \langle x, v_1\rangle > m - \langle z, v_1\rangle, |\langle x, v_j\rangle + \langle z, v_j\rangle| \le n \ \forall j \ne 1)].$$

It is not hard to check that this expression is asymptotic to a negative multiple of

$$n^{d-1} \sum_z p(z)\langle z, v_1\rangle.$$

Since the mean vector of $p(\cdot)$ is orthogonal to $v_1$, this is $o(n^{d-1})$, as required. □

REMARK. If one tried to prove Theorem 6 for general $v \in R^d$, one would need the continuity assumption in the definition of $v$-homogeneity for the coupling measure $\nu$, in order to know that, if $\nu$ concentrates on configurations with discrepancies of both types, then the rate of loss of discrepancies inside $B_{m,n}$ is of order $n^{d-1}$. As explained earlier, we do not know whether the continuity assumption on the marginals $\mu_i$ is inherited by the coupling measure $\nu$.

As a consequence of Theorem 6, we can explicitly identify all stationary $v$-profile measures when (6) is satisfied.

COROLLARY 1. *Suppose that* (6) *holds, and let*

$$\alpha_c(x) = \frac{ce^{\langle z, v\rangle}}{1 + ce^{\langle z, v\rangle}}, \qquad c > 0.$$

*Then, the set of all extremal stationary $v$-profile measures is exactly $\{\nu_{\alpha_c}, c > 0\}$.*

PROOF. Let $\mu$ be any extremal stationary $v$-profile measure. The measure $\nu_{\alpha_c}$ has the same properties: it is stationary by Theorem 2 and it is extremal by Theorem 2.1 of Jung [7], since

$$\sum_x \alpha_c(x)[1-\alpha_c(x)] = \infty.$$

(Jung's context assumes reversibility, but this property is not needed in his proof.) Therefore, by Theorem 6, for each $c > 0$, either $\mu \leq \nu_{\alpha_c}$ or $\nu_{\alpha_c} \leq \mu$. Since the measures $\nu_{\alpha_c}$ are stochastically monotone in $c$, there is some $0 \leq c^* \leq \infty$ so that $\mu \leq \nu_{\alpha_c}$ for $c > c^*$ and $\mu \geq \nu_{\alpha_c}$ for $c < c^*$. Since the family $\nu_{\alpha_c}$ is weakly continuous in $c$, it follows that $\mu = \nu_{\alpha_{c^*}}$. □

COROLLARY 2. *Suppose that*

$$\left\langle \sum_x x p(x), v \right\rangle = 0.$$

*Then, the extremal stationary $v$-homogeneous measures are exactly the homogeneous product measures.*

PROOF. The proof is the same as that of Corollary 1. The role that $\nu_{\alpha_c}$ played in the proof of Corollary 1 is now played by the homogeneous product measures. In this case, Theorem 1.17 on page 216 in [14] can be used as an alternative to Jung's theorem. □

Finally, we will combine Theorem 6 with part of the proof of Theorem 5 to show that, under an additional assumption on $p(\cdot)$, all $v$-profile stationary measures have a property analogous to that of a blocking measure in the sense of Bramson, Liggett and Mountford [1]. The proof is based on the proof of Lemma 6.4 of that paper. The simplest way to explain the idea of the proof of Theorem 7 is to refer to identity (6.12) in the statement of that lemma. It states that if $d = 1$ and $\nu$ is a stationary measure for the exclusion process with transition probabilities $p(\cdot)$ that has a finite expected number of ones to the left of the origin and a finite expected number of zeros to the right of the origin, then

$$\sum_{n=1}^{\infty} n^2 p(-n) = \sum_{n=1}^{\infty} n M(n)[p(n) - p(-n)],$$

where $M(n)$ is the $\nu$-expected number of sites $x$ at which there is a one at $x$ and a zero at $x + n$. This provides an a priori bound for the quantitites $M(n)$ in terms of the second moments of the negative jump probabilities for the particles. The proof below develops this idea for systems in higher dimensions. Most of the work is aimed at proving the convergence in (39).

THEOREM 7. *Suppose that*

(a) $p(z) \geq p(-z)$ *whenever* $\langle z, v \rangle > 0$,
(b) $p(z) > p(-z)$ *for some $z$ such that* $\langle z, v \rangle > 0$,

*and*

(c) $\sum_{z:\langle z,v \rangle > 0} p(-z) \langle z, v \rangle^2 < \infty$.

*If $\mu$ is a $v$-profile extremal stationary measure for the corresponding exclusion process, then $\mu$ is concentrated on configurations $\eta$ that satisfy $\eta(kv) = 0$ for all sufficiently large negative $k$ and $\eta(kv) = 1$ for all sufficiently large positive $k$.*

PROOF. We begin by restating (35), which only requires the assumption of stationarity. In the present case, $\Gamma = NZ$ for some integer $N \geq 1$, since $v \in Z^d$. We therefore have

$$\sum_z p(z)[f(kN) - f_z(kN)] = \sum_z p(z)[f(kN - \langle z, v \rangle) - f_z(kN - \langle z, v \rangle)].$$

Sum this identity over $m < k \leq n$. For fixed $z$, make the change of variables

$$k \to k + \frac{\langle z, v \rangle}{N}$$

in the sum over $k$ of the terms on the right-hand above. A significant amount of cancellation occurs with terms on the left-hand side, whose precise nature depends on the sign of $\langle z, v \rangle$. The result is

$$\sum_{z:\langle z,v \rangle > 0} p(z) \sum_{k=n+1-(\langle z,v \rangle/N)}^{n} [f(kN) - f_z(kN)]$$

$$- \sum_{z:\langle z,v \rangle < 0} p(z) \sum_{k=n+1}^{n-(\langle z,v \rangle/N)} [f(kN) - f_z(kN)]$$

$$= \sum_{z:\langle z,v \rangle > 0} p(z) \sum_{k=m+1-(\langle z,v \rangle/N)}^{m} [f(kN) - f_z(kN)]$$

$$- \sum_{z:\langle z,v \rangle < 0} p(z) \sum_{k=m+1}^{m-(\langle z,v \rangle/N)} [f(kN) - f_z(kN)].$$

Since the left-hand side above becomes the right-hand side if $n$ is replaced by $m$, it follows that the left-hand side is independent of $n$. Since $\mu$ is $v$-profile, this expression tends to zero as $n \to \pm\infty$. Therefore, it is zero for all $n$. In other words,

$$\sum_{z:\langle z,v \rangle > 0} p(z) \sum_{k=n+1-(\langle z,v \rangle/N)}^{n} [f(kN) - f_z(kN)]$$

$$= \sum_{z\,:\,\langle z,v\rangle<0} p(z) \sum_{k=n+1}^{n-(\langle z,v\rangle/N)} [f(kN) - f_z(kN)]$$

over all $n$. Next, sum this identity for $-M \leq n < M$. The result is

$$\sum_{z\,:\,\langle z,v\rangle>0} p(z) \sum_k [f(kN) - f_z(kN)]\left[\left(k + \frac{\langle z,v\rangle}{N}\right) \wedge M - (-M) \vee k\right]^+$$

$$= \sum_{z\,:\,\langle z,v\rangle<0} p(z) \sum_k [f(kN) - f_z(kN)]\left[k \wedge M - \left(k + \frac{\langle z,v\rangle}{N}\right) \vee (-M)\right]^+.$$

Replace $z$ by $-z$ on the right-hand above, then make a change of variables $k \to k + (\langle z,v\rangle/N)$, and finally use (36). The right-hand side above then becomes

$$\sum_{z\,:\,\langle z,v\rangle>0} p(-z) \sum_k [f(kN + \langle z,v\rangle) - f_z(kN)]\left[\left(k + \frac{\langle z,v\rangle}{N}\right) \wedge M - (-M) \vee k\right]^+.$$

Therefore, moving some terms from the resulting right-hand side to the left-hand side, we obtain

(38)
$$\sum_{z\,:\,\langle z,v\rangle>0} [p(z) - p(-z)] \sum_k [f(kN) - f_z(kN)] c_M(k,z)$$
$$= \sum_{z\,:\,\langle z,v\rangle>0} p(-z) \sum_k [f(kN + \langle z,v\rangle) - f(kN)] c_M(k,z),$$

where

$$c_M(k,z) = \left[\left(k + \frac{\langle z,v\rangle}{N}\right) \wedge M - (-M) \vee k\right]^+.$$

By comparing $\mu$ with a translate of $\mu$ in the $v$ direction, we see from Theorem 6 that $f(kN)$ is a monotone function of $k$. Since $\mu$ is a $v$-profile measure, it must be an increasing function of $k$. Therefore, all the summands on the right-hand side of (38) are nonnegative. Since

$$c_M(k,z) \leq \frac{\langle z,v\rangle}{N}$$

for all choices of the arguments, the right-hand side of (38) is at most

$$N^{-1} \sum_{z\,:\,\langle z,v\rangle>0} p(-z)\langle z,v\rangle \sum_k [f(kN + \langle z,v\rangle) - f(kN)]$$
$$= N^{-2} \sum_{z\,:\,\langle z,v\rangle>0} p(-z)\langle z,v\rangle^2,$$

which is finite by assumption. In the identity above, we have used the fact that $\langle z, v \rangle$ is a multiple of $N$.

The summands on the left-hand side of (38) are also nonnegative by assumption. Since

$$\lim_{M \to \infty} c_M(k, z) = \frac{\langle z, v \rangle}{N},$$

we may apply Fatou's lemma to the left-hand side of (38) to conclude that

$$\sum_{z : \langle z, v \rangle > 0} [p(z) - p(-z)] \langle z, v \rangle \sum_k [f(kN) - f_z(kN)] < \infty.$$

By assumption, there exists a $z$ with $\langle z, v \rangle > 0$ so that $p(z) > p(-z)$. For this $z$, it follows from this display that

(39) $$\sum_k [f(kN) - f_z(kN)] < \infty.$$

We claim that the convergence of the series in (39) does not depend on $z$, and hence it will converge for all $z$. Note first that by the definitions of $f$ and $f_z$,

(40) $$\mu\{\eta : \eta(x) = 1, \eta(x+z) = 0\} = f(\langle x, v \rangle) - f_z(\langle x, v \rangle).$$

To see that the convergence of (39) holds for all $z$, it suffices by (40) to check that, for any $z$ and $w$, there exists an $\varepsilon > 0$ so that, for all $x$,

$$\mu\{\eta : \eta(x) = 1, \eta(x+z) = 0\} \geq \varepsilon \mu\{\eta : \eta(x) = 1, \eta(x+w) = 0\}.$$

For this, it is enough to show that such an $\varepsilon$ exists so that

$$\mu\{\eta : \eta(x) = 1, \eta(x+z) = 0, \eta(x+w) = 1\}$$
$$\geq \varepsilon \mu\{\eta : \eta(x) = 1, \eta(x+w) = 0, \eta(x+z) = 1\}$$

for all $x$. To check this, since $\mu$ is stationary, it is enough to use irreducibility to find a path $z_0, z_1, \ldots, z_n$ with $z_0 = z$, $z_n = w$, and $p(z_{i+1} - z_i) > 0$ for each $i$, and then observe that there is a way to go from any configuration $\eta$ such that $\eta(x+z) = 1, \eta(x+w) = 0$ to $\eta_{x+z, x+w}$ (the configuration obtained from $\eta$ by interchanging the states at $x+z$ and $x+w$) by successively moving particles along the path without interference from sites not on the path (shifted by $x$), so that the probability is positive that these transitions occur in the given order by time one, and that by time one, no other transitions involving sites on this path are attempted. Schematically, suppose $\eta$ has the following form along the shifted path:

$$1 \quad 1 \quad 0 \quad 1 \quad 1 \quad 0 \quad 0$$

and we need to move to

$$0 \quad 1 \quad 0 \quad 1 \quad 1 \quad 0 \quad 1.$$

The way to do this is to move the rightmost 1 two steps to the right, the next 1 one step to the right, the next 1 two steps to the right, and finally the leftmost 1 one step to the right.

We have shown that the series in (39) converges for all $z$ and, in particular, for $z = v$. Recalling that $|v|^2$ is a multiple of $N$, say, $|v|^2 = lN$, we conclude that

$$\sum_k \mu\{\eta : \eta(kv) = 1, \eta((k+1)v) = 0\} = \sum_k [f(klN) - f_v(klN)] < \infty.$$

Therefore,

$$\lim_{k \to -\infty} \eta(kv) \quad \text{and} \quad \lim_{k \to +\infty} \eta(kv)$$

exist a.s. Since $\mu$ is a $v$-profile measure, these limits must be 0 and 1, respectively, as claimed. $\square$

**7. Convergence to stationary profile measures.** In previous sections we studied the stationary measures of exclusion processes. Here and in the next two sections, we give convergence results for exclusion processes with nonstationary initial measures. The exclusion processes are assumed to be on $Z^d$, $d \geq 2$, with transition probabilities that are translation invariant, have finite mean and are irreducible. The main result in this section is Theorem 8, where we show that, given a one-parameter family of extremal stationary $v$-profile measures, an exclusion process with a $v$-profile initial measure converges weakly to an average of these stationary measures, provided one is given certain bounds on the tails of the measures. As in Section 6, we assume $v$ is a nonzero element of $Z^d$.

In order to show Theorem 8, we will employ four lemmas. For the first of these, we introduce the following notation. Let $v_1, \ldots, v_d$ be an orthogonal basis for $R^d$ with $v_i \in Z^d$ and $v_1 = v$; as mentioned before Lemma 1, such a basis can be constructed using the Gram–Schmidt procedure. We set

$$L(z) = \{x \in Z^d : 0 \leq \langle x - z, v_i \rangle < \langle v_i, v_i \rangle, \ 2 \leq i \leq d\}$$

and

$$L^-(z) = L(z) \cap \{x \in Z^d : \langle x - z, v_1 \rangle \leq 0\},$$
$$L^+(z) = L(z) \cap \{x \in Z^d : \langle x - z, v_1 \rangle \geq 0\}.$$

That is, $L(z)$ is an infinite strip in the $v_1$ direction which has width $|v_i|$ in the perpendicular directions, and $L^-(z)$ and $L^+(z)$ are the corresponding semi-infinite strips. For $k = (k^{(2)}, \ldots, k^{(d)}) \in Z^{d-1}$, set

$$L_k(z) = L\left(z + \sum_{i=2}^d k^{(i)} v_i\right).$$

It is easy to see for given $z$ that $L_k(z)$, $k \in Z^{d-1}$, partitions $Z^d$ and

$$
(41) \qquad x \in L(z),\ y \in L_k(z) \implies \begin{cases} x - \sum_{i=2}^{d} k^{(i)} v_i \in L_{-k}(z), \\ y - \sum_{i=2}^{d} k^{(i)} v_i \in L(z). \end{cases}
$$

For a $v$-homogeneous measure $\lambda$, we also set

$$\ell_\lambda^-(z) = E^\lambda[\#x \in L^-(z) : \eta(x) = 1],$$
$$\ell_\lambda^+(z) = E^\lambda[\#x \in L^+(z) : \eta(x) = 0].$$

In the following lemmas and in Theorem 8, we will employ $v$-homogeneous measures $\lambda$ satisfying

$$(42\text{a}) \qquad \ell_\lambda^-(0) < \infty$$

and

$$(42\text{b}) \qquad \ell_\lambda^+(0) < \infty.$$

Note that such $v$-homogeneous measures are $v$-profile.

We will repeatedly use the coupling $\xi_t = (\eta_t, \zeta_t)$ introduced in Section 5 for the exclusion process. Unless stated otherwise, we will assume that the initial measures of the two coordinates are independent. Let $\mu$ and $\lambda$ denote the initial measures of the coordinates, respectively, and $\gamma$ the corresponding product measure.

Let $\xi_t$ be such a coupling, where $\mu$ and $\lambda$ are $v$-homogeneous measures satisfying (42b) for $\mu$ and (42a) for $\lambda$. The first part of Lemma 3 states that the expected number of discrepancies of type $(0,1)$, along a strip $L(0)$, is nonincreasing in time. The second part of the lemma uses this to obtain uniform upper bounds in time on the tail of $\ell_{\lambda_t}^-(z)$ for a measure $\lambda$ satisfying (42a), when a $v$-homogeneous stationary measure satisfying both parts of (42) exists.

LEMMA 3. *Assume that $\mu$ and $\lambda$ are $v$-homogeneous measures satisfying* (42b) *and* (42a), *respectively. Then,*

$$(43) \qquad E^\gamma[\#x \in L(0) : \eta_t(x) = 0,\ \zeta_t(x) = 1] \text{ is nonincreasing in } t.$$

*Assume, moreover, that $\mu$ is stationary and satisfies* (42a). *Then,*

$$(44) \qquad \ell_{\lambda_t}^-(z) \to 0 \qquad as\ \langle z, v \rangle \to -\infty$$

*uniformly in $t$.*

PROOF. Discrepancies of type $(0,1)$ can disappear, but cannot be created. So, in order to show (43), it suffices to show that the expected rate at which discrepancies leave $L(0)$ equals the expected rate at which they enter at a given time, if one ignores terms corresponding to their disappearance. Let $r(x, y; \xi)$ denote the rate at which a discrepancy of type $(0,1)$ moves from $x$ to $y$; this depends on just $p(y - x)$ and the present state $\xi$. The expected rate at which such discrepancies leave $L(0)$ at time $t$ is

$$E^\gamma\left[\sum_{x\in L(0)}\sum_{y\notin L(0)} r(x,y;\xi_t)\right] = \sum_{k\neq 0} E^\gamma\left[\sum_{x\in L(0)}\sum_{y\in L_k(0)} r(x,y;\xi_t)\right],$$

which is finite because of assumptions (42b) and (42a) on $\mu$ and $\lambda$. Since $\mu$ and $\lambda$ are $v$-homogeneous, so is $\gamma$, and hence so is the distribution of $\xi_t$. It therefore follows from (41) that the above sum equals

$$\sum_{k\neq 0} E^\gamma\left[\sum_{x\in L_{-k}(0)}\sum_{y\in L(0)} r(x,y;\xi_t)\right] = E^\gamma\left[\sum_{x\notin L(0)}\sum_{y\in L(0)} r(x,y;\xi_t)\right],$$

which is the expected rate at which these discrepancies enter $L(0)$.

We now demonstrate (44). Let $_n\mu$ denote the shift of $\mu$ by $-nv$ and $_n\gamma$ the product of $_n\mu$ and $\lambda$; since $\mu$ is $v$-homogeneous and satisfies (42b), the same is true for $_n\mu$. It follows from this and (43) applied to $_n\mu$ that, for given $\varepsilon > 0$ and sufficiently large $n$,

$$E^{_n\gamma}[\#x \in L(0): \eta_t(x) = 0,\ \zeta_t(x) = 1] < \varepsilon$$

uniformly in $t$.

On the other hand, since $\mu$ is stationary and is assumed to satisfy (42a),

$$E^{_n\mu}[\#x \in L^-(-2nv): \eta_t(x) = 1] = E^\mu[\#x \in L^-(-nv): \eta_t(x) = 1] < \varepsilon$$

for sufficiently large $n$ and all $t$. Together with the previous inequality, this implies

$$\ell^-_{\lambda_t}(-2nv) = E^\lambda[\#x \in L^-(-2nv): \zeta_t(x) = 1] < 2\varepsilon$$

for all $t$, since $L^-(-2nv) \subset L(0)$. This implies (44) for $z = -2nv$. The limit also holds for general $z$ since $L^-(z)$ is contained in a finite union of appropriate $L^-_{k_j}(-2nv)$, $j = 1, \ldots, 2^{d-1}$, where the individual coordinates of $k_j$ differ by at most 1 and $n$ is chosen so that

$$-2(n+1)\langle v, v\rangle \leq \langle z, v\rangle \leq -2n\langle v, v\rangle,$$

and since for any choice of $k$,

$$\ell^-_{\lambda_t}\left(-2nv + \sum_{i=2}^d k^{(i)} v_i\right) = \ell^-_{\lambda_t}(-2nv). \qquad \square$$

Let
$$L^N(0) = L(0) \cap \{x \in Z^d : -N\langle v,v \rangle \le \langle x,v \rangle < N\langle v,v \rangle\},$$
where $N \in Z^+$, and let $\hat{L}^N(0) = L(0) \setminus L^N(0)$. Also, let
$$B_N = \{x \in Z^d : -N\langle v_i, v_i \rangle \le \langle x, v_i \rangle < N\langle v_i, v_i \rangle \text{ for } i = 1, \ldots, d\},$$
$N \in Z^+$, where $\{v_1, \ldots, v_d\}$ is the orthogonal basis for $R^d$ defined earlier, with $v_1 = v$. Lemma 4 provides upper bounds on the number of discrepancies for the coupled process $\xi_t = (\eta_t, \zeta_t)$, starting from $v$-profile measures $\mu$ and $\lambda$ satisfying (42), where $\mu$ is stationary. The first part gives bounds on discrepancies of types $(0,1)$ and $(1,0)$ existing simultaneously in $B_N$, and the second part gives bounds on the total number of discrepancies in $\hat{L}^N(0)$.

LEMMA 4. *Assume that $\mu$ and $\lambda$ are $v$-profile measures satisfying* (42), *and that $\mu$ is stationary. Then,*

(45)   $P^\gamma\{\text{discrepancies of opposite types exist for } (\eta_t, \zeta_t) \text{ in } B_N\} \to 0$

*as $t \to \infty$ for each $N$. Also,*

(46)   $E^\gamma[\#\text{discrepancies of } (\eta_t, \zeta_t) \text{ in } \hat{L}^N(0)] \to 0 \quad \text{as } N \to \infty$

*uniformly in $t$.*

PROOF. In (45), we may assume that $N$ is large. Let $L_{N,B}$ denote the union of the $(4N)^{d-1}$ translations of $L(0)$ containing points in $B_{2N}$. It suffices to show for each $t$ at which the probability in (45) is at least $\varepsilon > 0$,

$E^\gamma[\#\text{discrepancies of } \xi_t \text{ in } L_{N,B}] - E^\gamma[\#\text{discrepancies of } \xi_{t+1} \text{ in } L_{N,B}] < -\varepsilon',$

where $\varepsilon' > 0$ depends on $\varepsilon$ and $N$, since this implies there are only a finite number of such $t$ at least distance 1 apart.

As in the demonstration of (43), one can show that the expected number of discrepancies in $L_{N,B}$ is nonincreasing at each $t$ by calculating the terms corresponding to their movement, but ignoring those corresponding to their disappearance. It therefore suffices to show that, for any configuration with a pair of discrepancies of opposite types at sites $z, w \in B_N$, there is a uniform positive lower bound on the probability of these discrepancies meeting each other or other discrepancies of opposite type after an additional unit of time, while remaining in $B_{2N} \subset L_{N,B}$. This is shown in Lemma 3.1 in [1]. The argument uses the irreducibility of $p(\cdot)$, and is similar to that near the end of the proof of Theorem 7, where one constructs a path $z_0, z_1, \ldots, z_n$, with $z_0 = z$, $z_n = w$, and $p(z_{i+1} - z_i) > 0$ for each $i$, along which the discrepancies may move while remaining in $B_{2N}$, until meeting one another or some other discrepancy of opposite type.

To show (46), we apply (44) both to $\lambda$ and to the corresponding measure with the roles of 0 and 1 reversed, to get

$$\ell^-_{\lambda_t}(-z) \to 0 \quad \text{and} \quad \ell^+_{\lambda_t}(z) \to 0 \qquad \text{as } \langle z, v \rangle \to \infty$$

uniformly in $t$. Since $\mu$ is stationary and is assumed to satisfy (42), $\mu_t$ satisfies the analogous limits. Together, these four limits imply (46). $\square$

The following elementary lemma compares extremal stationary measures for the exclusion process on the sets

$$D_M = \{(\eta, \zeta) : \eta(z) \geq \zeta(z) \text{ for } z \in B_M\},$$

where $M \in Z^+$.

LEMMA 5. *Let $\gamma$ be any measure with marginals $\mu^{\alpha_1}$ and $\mu^{\alpha_2}$, where $\mu^{\alpha_1}$ and $\mu^{\alpha_2}$ are extremal stationary measures with $\mu^{\alpha_1}\{\eta : \eta(0) = 1\} < \mu^{\alpha_2}\{\eta : \eta(0) = 1\}$. Then, for any $\varepsilon > 0$,*

(47) $$\gamma(D_M) < \varepsilon$$

*for large enough $M$ not depending on the choice of $\gamma$.*

PROOF. Assume, on the contrary, that there exist such measures $\gamma^k$ and sets $D_{M_k}$ with $\gamma^k(D_{M_k}) \geq \varepsilon$ and $M_k \to \infty$, as $k \to \infty$. Choose a subsequence $k_n$ along which $\gamma^{k_n} \xrightarrow{w} \gamma^\infty$ as $n \to \infty$ for some measure $\gamma^\infty$; then $\gamma^\infty$ has the same marginals $\mu^{\alpha_i}$. Let

$$D = \{(\eta, \zeta) : \eta(z) \geq \zeta(z) \text{ for all } z\};$$

we claim that $\gamma^\infty(D) = 0$. It will follow that $\gamma^{k_n}(D_M) < \varepsilon$ for large enough $n$ and $M$, which will contradict our assumption on $\gamma^k$.

To see that $\gamma^\infty(D) = 0$, consider the coupled processes $\eta_t$ and $\zeta_t$ with initial joint measure $\gamma^\infty$. Let $\bar{\gamma}_t$ denote the average of $\gamma_s^\infty$ over $s \in [0, t]$, and $\bar{\gamma}$ the weak limit of $\bar{\gamma}_{t_n}$ along some subsequence $t_n$. Then, $\bar{\gamma}$ is stationary with marginals $\mu^{\alpha_i}$, and $D$ is invariant in time under $\bar{\gamma}$. Since $\mu^{\alpha_i}$ are assumed to be extremal with densities $\alpha_1 < \alpha_2$ at 0, one has $\bar{\gamma}(D) = 0$. Because $\gamma^\infty(D) \leq \bar{\gamma}(D)$, one has $\gamma^\infty(D) = 0$, as desired. $\square$

In the following lemma and Theorem 8, we will assume that

(48) there exists a one-parameter family $\mu^\alpha$, $\alpha \in (0, 1)$, of extremal stationary $v$-profile measures with $\mu^\alpha\{\eta : \eta(0) = 1\} = \alpha$.

It follows immediately from part (a) of Theorem 6 that

$\mu^\alpha, \alpha \in (0, 1)$, are stochastically ordered, with $\mu^{\alpha_1} \leq \mu^{\alpha_2}$ for $\alpha_1 \leq \alpha_2$.

For $\alpha_0 \in (0,1)$ and $z \in Z^d$, one also has

(49) $\quad\quad\quad \mu^\alpha\{\eta : \eta(z) = 1\} \to \mu^{\alpha_0}\{\eta : \eta(z) = 1\} \quad\quad \text{as } \alpha \to \alpha_0.$

This follows by noting that, as either $\alpha \nearrow \alpha_0$ or $\alpha \searrow \alpha_0$, $\mu^\alpha \xrightarrow{w} \mu'$, where $\mu'$ is stationary with $\mu'\{\eta : \eta(0) = 1\} = \alpha_0$. When $\alpha \nearrow \alpha_0$ $(\alpha \searrow \alpha_0)$, one has $\mu' \leq \mu^{\alpha_0}$ $(\mu' \geq \mu^{\alpha_0})$; elementary reasoning similar to the path argument near the end of the proof of Theorem 7 will then imply $\mu' = \mu^{\alpha_0}$, and hence (49). (Otherwise, the coupling measure for $\mu'$ and $\mu^{\alpha_0}$ would have discrepancies at some sites, but not at 0.) From Theorem 6 and (49), it follows that $\mu^\alpha$, $\alpha \in (0,1)$, are the unique extremal stationary $v$-profile measures. One can also check that, for each $z \in Z^d$,

(50) $\quad\quad\quad \begin{aligned} \mu^\alpha\{\eta : \eta(z) = 1\} &\to 0 \quad\quad \text{as } \alpha \to 0, \\ &\to 1 \quad\quad \text{as } \alpha \to 1, \end{aligned}$

since this holds at $z = 0$.

Like Lemma 5, the following lemma compares stationary measures. We consider here a stationary coupling of $\mu^{\alpha_0}$ with $\int_{(0,1)} \mu^\alpha \sigma(d\alpha)$, where $\mu^{\alpha_0}$ and $\mu^\alpha$ satisfy (48), and the configurations almost everywhere have discrepancies of at most one type. Then, $\mu^{\alpha_0}$ will "dominate" and "be dominated by" $\int_{(0,1)} \mu^\alpha \sigma(d\alpha)$ on invariant subsets corresponding to $\alpha \in (0, \alpha_0]$ and $\alpha \in (\alpha_0, 1)$, respectively. We employ the notation

$$D^- = \{(\eta, \zeta) : \eta(z) < \zeta(z) \text{ for some } z \in Z^d\},$$
$$D^+ = \{(\eta, \zeta) : \eta(z) > \zeta(z) \text{ for some } z \in Z^d\}.$$

LEMMA 6. *Let the coupling $\xi_t = (\eta_t, \zeta_t)$ be stationary, with measure $\gamma$ having marginals $\mu^{\alpha_0}$, $\alpha_0 \in (0,1)$, and $\int_{(0,1)} \mu^\alpha \sigma(d\alpha)$, where $\mu^{\alpha_0}$ and $\mu^\alpha$ are $v$-profile measures satisfying* (48), *and $\sigma$ is a probability measure. Also, let $\gamma$ be concentrated on configurations with discrepancies of at most one type. Then, there exist invariant subsets $G^- = \{0,1\}^{Z^d} \times H^-$ and $G^+ = \{0,1\}^{Z^d} \times H^+$ that partition the space, with $\int_{(0,\alpha_0]} \mu^\alpha \sigma(d\alpha)$ and $\int_{(\alpha_0,1)} \mu^\alpha \sigma(d\alpha)$ concentrated on $H^-$ and $H^+$, respectively. Moreover,*

$$\gamma(D^- \cap G^-) = \gamma(D^+ \cap G^+) = 0.$$

PROOF. The existence of invariant sets $H^-$ and $H^+$ that partition $\{0,1\}^{Z^d}$ and are concentrated on $\int_{(0,\alpha_0]} \mu^\alpha \sigma(d\alpha)$ and $\int_{(\alpha_0,1)} \mu^\alpha \sigma(d\alpha)$ follows, with a little work, from the unique representation of the stationary measures in terms of their extremal elements. It is easy to see that $G^\pm = \{0,1\}^{Z^d} \times H^\pm$ partition $\{0,1\}^{Z^d} \times \{0,1\}^{Z^d}$ and are invariant on $\gamma$.

We will show that $\gamma(D^- \cap G^-) = 0$; the argument for $\gamma(D^+ \cap G^+)$ is analogous. Assume, on the contrary, that $\gamma(D^- \cap G^-) > 0$, and denote by $\gamma^-$, the stationary measure obtained by conditioning on $D^- \cap G^-$. By the extremality of $\mu^{\alpha_0}$, the first marginal measure of $\gamma^-$ is $\mu^- = \mu^{\alpha_0}$. Also, by the extremality of $\mu^\alpha$, $\alpha \in (0, \alpha_0]$, the second marginal measure $\lambda^-$ can be written as

$$\lambda^- = \int_{(0,\alpha_0]} \mu^\alpha \sigma^-(d\alpha)$$

for some probability measure $\sigma^-$. This implies

$$\mu^-\{\eta : \eta(z) = 1\} \geq \lambda^-\{\zeta : \zeta(z) = 1\}$$

for all $z$. However, since $\gamma^-$ is supported on $D^-$, this inequality cannot hold for all $z$, which produces a contradiction. Hence, $\gamma(D^- \cap G^-) = 0$, as desired. □

In the proof of Theorem 8, we will employ an extension $\vec{\xi}_t = (\vec{\eta}_t, \zeta_t)$ of the coupling $\xi_t = (\eta_t, \zeta_t)$ introduced already, where the coordinates of $\vec{\eta}_t = (\eta_t^\alpha)_{\alpha \in (0,1)}$ have initial measures $\mu^\alpha$ which satisfy (48), with $\eta_0^{\alpha_1}(z) \leq \eta_0^{\alpha_2}(z)$ for $\alpha_1 \leq \alpha_2$ and all $z$. It then follows that $\eta_t^{\alpha_1}(z) \leq \eta_t^{\alpha_2}(z)$ for all $t$. We denote by $\lambda$ the measure of $\zeta_0$, by $\vec{\mu}$ the measure of $\vec{\eta}_0$ and by $\vec{\gamma}$ the measure of $\vec{\xi}_0 = (\vec{\eta}_0, \zeta_0)$. As before, we assume $\vec{\eta}_0$ and $\zeta_0$ are independent.

Theorem 8 is the main result of the section. It states that when there exist $v$-profile stationary measures $\mu^\alpha$ satisfying (42) and (48), and a $v$-profile measure $\lambda$ satisfying (42), then the measure of the process starting at $\lambda$ will converge to an average of the former as $t$ goes to infinity. In one dimension, this statement is an immediate consequence of the convergence theorem for positive recurrent Markov chains. In higher dimensions, quite a bit more work is required to prove it.

THEOREM 8. *Assume that there exist $v$-profile measures $\mu^\alpha$, $\alpha \in (0,1)$, that are stationary for the exclusion process and satisfy (42) and (48). Also, assume that (42) is satisfied by a $v$-profile measure $\lambda$. Then,*

(51) $$\lambda_t \overset{w}{\to} \lambda_\infty \qquad as\ t \to \infty$$

*for some stationary $v$-profile measure $\lambda_\infty$, where*

(52) $$\lambda_\infty = \int_{(0,1)} \mu^\alpha \sigma(d\alpha)$$

*and $\sigma$ is a probability measure.*

PROOF. The argument consists of two main parts. We will first show that any limit of $\lambda_{t_k}$, along a subsequence $t_k$, must be of the form (52). We will then show that $\sigma$ in (52) does not depend on the subsequence, and so this limit must in fact hold along the entire sequence. The first part of the argument uses (44) of Lemma 3, (45) of Lemma 4 and Lemma 5; the second part uses (43) of Lemma 3, both parts of Lemma 4 and Lemma 6.

*Characterization of limits.* We will first show that any limit can be approximated locally by a mixture of $\mu^\alpha$, $\alpha = j/n$, $j = 1, \ldots, n-1$, if each $\mu^\alpha$ is restricted to an appropriate subset of configuration space. We will then show that conditioning on this restriction does not change the distribution of $\mu^\alpha$ much. Such a limit will therefore be close to an average of the unrestricted $\mu^\alpha$. Taking limits as $n \to \infty$ will imply (52).

We apply here the coupling that was introduced between $\eta_t^\alpha$, $\alpha \in (0,1)$, and $\zeta_t$, with initial measures $\mu^\alpha$ and $\lambda$. By (45) of Lemma 4, for each $N$ and choice of $\alpha = j/n$, $j = 1, \ldots, n-1$, $(\eta_t^\alpha, \zeta_t)$ typically does not have discrepancies of both types $(0,1)$ and $(1,0)$ on $B_N$, for large enough $t$. Also, on account of (44) of Lemma 3, its analog for $\ell_{\lambda_t}^+(z)$, and (50),

$$P^{\vec{\gamma}}\{\vec{\xi} : \eta_t^\alpha(z) \geq \zeta_t(z) \text{ for all } z \in L^N(0)\} \to 0 \qquad \text{as } \alpha \to 0,$$

$$P^{\vec{\gamma}}\{\vec{\xi} : \eta_t^\alpha(z) \leq \zeta_t(z) \text{ for all } z \in L^N(0)\} \to 0 \qquad \text{as } \alpha \to 1$$

and $N \to \infty$, uniformly in $t$; since $L^N(0) \subset B_N$, we may replace $L^N(0)$ with $B_N$ in the display. It follows that, for given $\varepsilon > 0$, large $n, N$ and small $a_\varepsilon > 0$,

$$P^{\vec{\gamma}}\left(\bigcup_{j \in A_{\varepsilon,n}} {}^N F_t^{j,n}\right) \geq 1 - \varepsilon \tag{53}$$

for large enough $t$, where $A_{\varepsilon,n} = \{j : j \in [a_\varepsilon n + 1, (1 - a_\varepsilon)n - 1]\}$ and where

$${}^N F_t^{j,n} = \{\vec{\xi} : \eta_t^{j,n}(z) \leq \zeta_t(z) \leq \eta_t^{j+1,n}(z) \text{ for all } z \in B_N\} \cap \left(\bigcup_{j'=1}^{j-1} {}^N F_t^{j',n}\right)^c.$$

That is, off of a set of probability $\varepsilon$, the space can be partitioned into sets ${}^N F_t^{j,n}$, $j \in A_{\varepsilon,n}$, such that, for each $j$, the configuration $\zeta_t$ is sandwiched between $\eta_t^{j,n}$ and $\eta_t^{j+1,n}$ on $B_N$. (Here and later on, when $\alpha = j/n$, we often write $\eta_t^{j,n}$ for $\eta^\alpha$.)

We do not know much about the behavior of $\eta_t^{j,n}$ on ${}^N F_t^{j,n}$. We will show, however, that $\eta_{t'}^{j,n}$ is almost independent of ${}^N F_t^{j,n}$ for $t' = t + s$, where $s$ is large but small relative to $N$, if we consider only $z \in B_M$, for fixed $M$. First note that, when $s \in [0, s_0]$ for given $s_0$, if $N$ is chosen large enough, it follows from (53) and the shift invariance of the transition probabilities that

$$P^{\vec{\gamma}}\left(\bigcup_{j \in A_{\varepsilon,n}} \{\vec{\xi} \in {}^N F_t^{j,n} : \eta_{t'}^{j,n}(z) \leq \zeta_{t'}(z) \leq \eta_{t'}^{j+1,n}(z) \text{ for all } z \in B_M\}\right)$$

(54)
$$\geq 1 - 2\varepsilon,$$

for given $M$. That is, on ${}^N F_t^{j,n}$, the configuration $\zeta_{t'}$ typically remains sandwiched between $\eta_{t'}^{j,n}$ and $\eta_{t'}^{j+1,n}$ on $B_M$. Let ${}^M \lambda_{t'}$ denote the measure obtained from $\zeta_{t'}$ by restricting configurations from $Z^d$ to $B_M$, let ${}^{M,N}\mu_{t',t}^{j,n}$ denote the measure obtained from $\eta_{t'}^{j,n}$ by also conditioning on ${}^N F_t^{j,n}$, and set ${}^N c_t^{j,n} = P^{\vec{\gamma}}({}^N F_t^{j,n})$. Also, let $d_M(\cdot, \cdot)$ denote the total variational distance between measures of configurations on $B_M$. (Since the space is finite, this is just the sum of the absolute value of the differences of the probabilities over all configurations.) Employing this notation, we can restate (54) as

$$\text{(55a)} \quad \sum_{j \in A_{\varepsilon,n}} {}^N c_t^{j,n} \, {}^{M,N}\mu_{t',t}^{j,n} \leq {}^{M,N}\tilde{\lambda}_{t'} \leq \sum_{j \in A_{\varepsilon,n}} {}^N c_t^{j,n} \, {}^{M,N}\mu_{t',t}^{j+1,n},$$

for some ${}^{M,N}\tilde{\lambda}_{t'}$ with

$$\text{(55b)} \quad d_M({}^M\lambda_{t'}, {}^{M,N}\tilde{\lambda}_{t'}) \leq 2\varepsilon.$$

On the other hand, for each $j$, one can consider the coupled process $(\eta_s', \zeta_s')$, where the initial measures are given by $\mu^{j-1,n}$ and $\mu^{j,n}$ and the coordinates are initially independent with joint measure $\gamma'$. As before, it follows from (45) of Lemma 4 that, for given $M$, $(\eta_s, \zeta_s)$ typically does not have discrepancies of both types on $B_M$ for large $s$. Moreover, since both measures are extremal stationary and the density of $\mu^{j,n}$ at 0 is greater than that of $\mu^{j-1,n}$, it follows from this and Lemma 5 that, for large $s$,

$$P^{\gamma'}\{(\eta', \zeta') : \eta_s'(z) \leq \zeta_s'(z) \text{ for all } z \in B_M\} \geq 1 - \varepsilon'$$

for given $\varepsilon' > 0$. Set $\zeta_0' = \eta_t^{j,n}$ and condition on the set ${}^N F_t^{j,n}$. For such $s$, it follows that, for $\varepsilon'$ small relative to ${}^N c_t^{j,n}$, $\eta_{t'}^{j,n}$ typically dominates, at sites in $B_M$, a process with measure $\mu^{j-1,n}$ (since conditioning on ${}^N F_t^{j,n}$ does not change the measure of the process starting at $\eta_0'$). The same reasoning also shows that, after conditioning on ${}^N F_t^{j,n}$, $\eta_{t'}^{j+1,n}$ is, for large $s$, typically dominated, on $B_M$, by a process with measure $\mu^{j+2,n}$, when $\varepsilon'$ is small relative to ${}^N c_t^{j,n}$. This domination translates into the inequalities

$$\text{(56a)} \quad {}^M\mu^{j-1,n} \leq {}^{M,N}\tilde{\mu}_{t',t}^{j,n} \quad \text{and} \quad {}^{M,N}\tilde{\mu}_{t',t}^{j+1,n} \leq {}^M\mu^{j+2,n}$$

for some ${}^{M,N}\tilde{\mu}_{t',t}^{j,n}$ and ${}^{M,N}\tilde{\mu}_{t',t}^{j+1,n}$ with

$$\text{(56b)} \quad {}^N c_t^{j,n} d_M({}^{M,N}\mu_{t',t}^{j,n}, {}^{M,N}\tilde{\mu}_{t',t}^{j,n}) \leq \varepsilon', \quad {}^N c_t^{j,n} d_M({}^{M,N}\mu_{t',t}^{j+1,n}, {}^{M,N}\tilde{\mu}_{t',t}^{j+1,n}) \leq \varepsilon'.$$

Setting $\varepsilon' = \varepsilon/n$, (55) and (56) together imply that, for given $\varepsilon > 0$ and large $n$, $N$,

$$\text{(57a)} \quad \sum_{j \in A_{\varepsilon,n}} {}^N c_t^{j,n} {}^M\mu^{j-1,n} \leq {}^{M,N}\tilde{\lambda}_{t'} \leq \sum_{j \in A_{\varepsilon,n}} {}^N c_t^{j,n} {}^M\mu^{j+2,n}$$

for some ${}^{M,N}\tilde{\lambda}_{t'}$ with

$$\text{(57b)} \quad d_M({}^M\lambda_{t'}, {}^{M,N}\tilde{\lambda}_{t'}) \leq 4\varepsilon,$$

if $s = t' - t$ is large relative to $M$, but small relative to $t'$. Suppose now that $\lambda_{t'_k} \overset{w}{\to} \lambda_\infty$ as $t'_k \to \infty$, for some measure $\lambda_\infty$. Then, the analogue of (57a) holds with $n = n_k \to \infty$ and $N = N_k \to \infty$, but with (57b) replaced by

$$d_M(^M\lambda_{t'_k}, {}^{M,N}\tilde{\lambda}_{t'_k}) \to 0 \qquad \text{as } k \to \infty.$$

Since $d_M(^M\mu^{j-1,n}, {}^M\mu^{j+2,n}) \to 0$ uniformly in $j$ as $n \to \infty$, it follows, with a little work, that

$$^M\lambda_\infty = \int_{(0,1)} {}^M\mu^\alpha \sigma(d\alpha),$$

where $\sigma$ is the weak limit of $\sum_{j \in A_{\varepsilon,n}} {}^N c_t^{j,n} \delta_{j/n}$ along some subsequence $t'_{k_i}$; because the measures $\sum_{j \in A_{\varepsilon,n}} {}^N c_t^{j,n} \delta_{(j-1)/n}$ and $\sum_{j \in A_{\varepsilon,n}} {}^N c_t^{j,n} \delta_{(j+2)/n}$ are concentrated on $[a_\varepsilon, 1 - a_\varepsilon]$, the sequence is tight on $(0,1)$. ($\delta_\alpha$ denotes the point mass at $\alpha$.) Letting $M \to \infty$, one obtains

$$\lambda_\infty = \int_{(0,1)} \mu^\alpha \sigma(d\alpha),$$

which is (52). It follows from this representation for $\lambda_\infty$ [or, alternatively, from (44)] that $\lambda_\infty$ is a $v$-profile measure.

*Uniqueness of $\sigma$.* We now show that weak limits $\lambda_\infty$ and $\lambda'_\infty$ of $\lambda_t$ along different subsequences must, in fact, have the same measures $\sigma$ and $\sigma'$ in (52). Let $g(\alpha_0, \alpha)$, with $\alpha_0, \alpha \in (0,1)$, denote the expected number of discrepancies of type $(0,1)$ on $L(0)$ for the pair $(\eta, \zeta)$ with marginal measures $\mu^{\alpha_0}$ and $\mu^\alpha$, where the joint measure is concentrated on configurations satisfying

$$
\begin{aligned}
\eta(z) &\geq \zeta(z) \qquad \text{for all } z, \text{ if } \alpha_0 \geq \alpha, \\
&\leq \zeta(z) \qquad \text{for all } z, \text{ if } \alpha_0 \leq \alpha.
\end{aligned}
\tag{58}
$$

The first part of the argument will be to show that

$$\int_{(0,1)} g(\alpha_0, \alpha) \sigma(d\alpha) = \int_{(0,1)} g(\alpha_0, \alpha) \sigma'(d\alpha) \tag{59}$$

for each $\alpha_0$.

We consider the coupled process $(\eta_\cdot, \zeta_\cdot)$ with marginal initial measures $\mu^{\alpha_0}$, $\alpha_0 \in (0,1)$, and $\lambda$. In addition to assuming $\lambda_{t_k} \overset{w}{\to} \lambda_\infty$ as $k \to \infty$, with $\lambda_\infty$ satisfying (52), we may also assume that the joint measures $\gamma_{t_k} \overset{w}{\to} \gamma_\infty$, for some $\gamma_\infty$. Let $\gamma_t^\infty$ denote the process restarted from $\gamma_0^\infty = \gamma_\infty$, and let $\bar{\gamma}_t$ be the Cesaro average of $\gamma_s^\infty$ over $s \in [0, t]$. Any weak limit $\bar{\gamma}$ of $\bar{\gamma}_t$, as $t \to \infty$, is stationary, with the same marginals $\mu^{\alpha_0}$ and $\int_{(0,1)} \mu^\alpha \sigma(d\alpha)$. Also, by (46) of Lemma 4, $\gamma_\infty$ and, hence, $\bar{\gamma}$, is concentrated on configurations with discrepancies of at most one type.

Let $g_t(\alpha_0)$ denote the expected number of discrepancies of type $(0,1)$ on $L(0)$ at time $t$, $g_t^N(\alpha_0)$ the corresponding expected number on $\hat{L}^N(0)$, and $\bar{g}(\alpha_0)$ the expected number on $L(0)$ for the measure $\bar{\gamma}$. By (43) of Lemma 3, $g_t(\alpha_0)$ is nonincreasing in $t$, so $g_\infty(\alpha_0) \stackrel{\text{def}}{=} \lim_{t\to\infty} g_t(\alpha_0)$ exists. Also, by (46) of Lemma 4,

$$g_t^N(\alpha_0) \to 0 \qquad \text{as } N \to \infty$$

uniformly in $t$. Therefore, since $\gamma_{t_k} \stackrel{w}{\to} \gamma_\infty$, it follows that

$$\bar{g}(\alpha_0) = g_\infty(\alpha_0).$$

We claim that

(60) $$\int_{(0,1)} g(\alpha_0, \alpha)\sigma(d\alpha) = \bar{g}(\alpha_0).$$

Since $g_\infty(\alpha_0)$ does not depend on the choice of $t_k$, it will follow from the previous two equations that

$$\int_{(0,1)} g(\alpha_0, \alpha)\sigma(d\alpha) = \int_{(0,1)} g(\alpha_0, \alpha)\sigma'(d\alpha) = g_\infty(\alpha_0),$$

and so (59) in fact holds.

To show (60), note that

$$\bar{g}(\alpha_0) = E^{\bar{\gamma}}[\#z \in L(0) : \eta(z) < \zeta(z); G^-] + E^{\bar{\gamma}}[\#z \in L(0) : \eta(z) < \zeta(z); G^+],$$

where we choose $G^-$ and $G^+$ as in Lemma 6. By Lemma 6, the first term on the right-hand side is 0. The second term on the right-hand side equals

$$\sum_{z \in L(0)} \left[ \int_{(\alpha_0, 1)} \mu^\alpha\{\eta : \eta(z) = 1\}\sigma(d\alpha) - \mu^{\alpha_0}\{\eta : \eta(z) = 1\}\sigma\{\alpha : \alpha > \alpha_0\} \right]$$

$$+ E^{\bar{\gamma}}[\#z \in L(0) : \eta(z) > \zeta(z); G^+].$$

By Lemma 6, the above expectation is 0. It follows that

$$\bar{g}(\alpha_0) = \sum_{z \in L(0)} \left[ \int_{(\alpha_0, 1)} \mu^\alpha\{\eta : \eta(z) = 1\}\sigma(d\alpha) - \mu^{\alpha_0}\{\eta : \eta(z) = 1\}\sigma\{\alpha : \alpha > \alpha_0\} \right].$$

On the other hand, it is easy to see that, because of the coupling in (58), $\int_{(\alpha_0, 1)} g(\alpha_0, \alpha)\sigma(d\alpha)$ is equal to the right-hand side of the last equation and

$$\int_{(0,\alpha_0]} g(\alpha_0, \alpha)\sigma(d\alpha) = 0.$$

Putting these terms together, we see that (60) holds, which implies (59).

We will show that (59) is not possible unless $\sigma = \sigma'$. We first note that, for $h \in (0, \alpha_0)$,

$$
\begin{aligned}
g(\alpha_0 - h, \alpha) - g(\alpha_0, \alpha) &= g(\alpha_0 - h, \alpha_0) && \text{for } \alpha \geq \alpha_0, \\
&= 0 && \text{for } \alpha < \alpha_0 - h.
\end{aligned}
\tag{61}
$$

The first equality holds since uncoupled particles from the first term on the left-hand side arise from either uncoupled particles from the second term on the left-hand side or from the term on the right-hand side; the second inequality holds trivially since both terms on the left-hand side are 0. Consequently, for each $\alpha_0 \in (0,1)$ and $h \in (0, \alpha_0)$,

$$\int_{[\alpha_0, 1)} (g(\alpha_0 - h, \alpha) - g(\alpha_0, \alpha)) \sigma(d\alpha) = g(\alpha_0 - h, \alpha_0) \sigma([\alpha_0, 1)),$$

with the analogous equality also holding for $\sigma'$. Also, note that for $\alpha \in [\alpha_0 - h, \alpha_0)$,

$$0 \leq g(\alpha_0 - h, \alpha) - g(\alpha_0, \alpha) \leq g(\alpha_0 - h, \alpha) \leq g(\alpha_0 - h, \alpha_0).$$

The last two displays, together with the second part of (61), imply that

$$\int_{(0,1)} (g(\alpha_0 - h, \alpha) - g(\alpha_0, \alpha)) \sigma(d\alpha) - \int_{(0,1)} (g(\alpha_0 - h, \alpha) - g(\alpha_0, \alpha)) \sigma'(d\alpha)$$
$$\leq g(\alpha_0 - h, \alpha_0)(\sigma([\alpha_0 - h, 1)) - \sigma'([\alpha_0, 1))).$$

By (59), the left-hand side of this inequality is 0, and so the right-hand side is nonnegative. Note that $g(\alpha_0 - h, \alpha_0)$ is at least the difference of the probabilities of there being a particle at 0 for $\mu^{\alpha_0}$ and for $\mu^{\alpha_0 - h}$, which by assumption is $h$. So,

$$\sigma([\alpha_0 - h, 1)) \geq \sigma'([\alpha_0, 1)).$$

Letting $h \searrow 0$ implies

$$\sigma([\alpha_0, 1)) \geq \sigma'([\alpha_0, 1)) \qquad \text{for each } \alpha_0 \in (0,1).$$

Since the reverse inequality also holds, one has $\sigma = \sigma'$, as desired. □

When (6) holds, the product measures $\nu_{\alpha_c}$, $c > 0$, in Corollary 1 of Theorem 6, satisfy the conditions in (48); it is easy to see that they also satisfy (42). One thus has the following immediate consequence of Theorem 8.

COROLLARY 3. *Assume that the kernel $p(\cdot)$ satisfies (6), and that a given v-profile measure $\lambda$ satisfies (42). Then,*

$$\lambda_t \xrightarrow{w} \lambda_\infty \qquad \text{as } t \to \infty$$

*for some stationary v-profile measure $\lambda_\infty$ satisfying (52).*

**8. Hydrodynamics and local limits.** Hydrodynamic scaling has been applied to a wide variety of stochastic processes, including the exclusion processes on $Z^d$, to obtain detailed information about their asymptotic behavior. In the setting of the exclusion processes, one obtains deterministic limits that satisfy the well-known Burgers' equation. These results extend to local limits for the unscaled exclusion processes, that are product measures with densities given by the solution of Burgers' equation there. After stating these results, we provide several concrete examples of such limits for exclusion processes in $d = 2$. We then apply hydrodynamic scaling to obtain a quick (though based on some nontrivial results) alternative proof of Theorem 4 (with the weak inequality). As in previous sections, the exclusion processes are assumed to be on $Z^d$, with transition probabilities that are translation invariant and irreducible; here, they are also assumed to have finite range. For the reader's convenience, certain basic results pertaining to entropy solutions of Burgers' equation are given in the Appendix.

The time evolution of the asymmetric exclusion processes on $Z$ with nearest transition probabilities $p(\cdot)$ has been extensively studied for product initial measures with constant densities $r$ and $\ell$ over the positive and negative half-lines. (See [14] for references.) Setting $p(1) = 1 - p(-1) = p$ and letting $\lambda_t$ denote the corresponding measure at time $t$, one has, for $p > 1/2$,

$$
(62) \qquad \lim_{t \to \infty} \lambda_t = \begin{cases} \nu_{1/2}, & \text{for } r \leq \tfrac{1}{2} \text{ and } \ell \geq \tfrac{1}{2}, \\ \nu_r, & \text{for } r \geq \tfrac{1}{2} \text{ and } \ell + r > 1, \\ \nu_\ell, & \text{for } \ell \leq \tfrac{1}{2} \text{ and } \ell + r < 1, \\ \tfrac{1}{2}\nu_\ell + \tfrac{1}{2}\nu_r, & \text{for } 0 < \ell < r \text{ and } \ell + r = 1. \end{cases}
$$

On the level of heuristics, the first three lines of (62) can be motivated by using approximations leading to the one-dimensional Burgers' equation

$$
(63a) \qquad \frac{\partial u}{\partial t} + m \frac{\partial}{\partial y}[u(1-u)] = 0,
$$

with

$$
(63b) \qquad \begin{aligned} u(0,y) &= r \quad \text{for } y \geq 0, \\ &= \ell \quad \text{for } y < 0, \end{aligned}
$$

where $m$ is the mean of $p(\cdot)$. Under this initial data, the entropy solution of (63) evolves in two different ways depending on whether or not $\ell < r$. If $\ell < r$, then $u(t, \cdot)$ is given by the *shock wave* that is a translate of $u(0, \cdot)$ in (63b), with

$$
(64) \qquad \begin{aligned} u(t,y) &= r \quad \text{for } y \geq m(1 - r - \ell)t, \\ &= \ell \quad \text{for } y < m(1 - r - \ell)t. \end{aligned}
$$

If $\ell > r$, then $u(t, y)$ is a *rarefaction wave* that is continuous for $t > 0$, with

$$u(t, y) = r \quad \text{for } y \geq m(1 - 2r)t,$$
$$= \ell \quad \text{for } y \leq m(1 - 2\ell)t, \tag{65}$$

and is linear over $[m(1 - 2\ell), m(1 - 2r)]$. Where $u(t, \cdot)$ is locally nearly constant, it is reasonable to expect $\lambda_t$ to be close to a product measure. Substitution of $y = 0$ into (64) and (65) then gives the densities in (62).

Hydrodynamic scaling provides a rigorous connection between entropy solutions of Burgers' equation and the asymptotic behavior of the exclusion processes in $d$ dimensions. For exclusion processes, the relevant formulation of Burgers' equation in $d$ dimensions is

$$\frac{\partial u}{\partial t} + \sum_{i=1}^{d} m_i \frac{\partial}{\partial x_i}[u(1-u)] = 0, \tag{66a}$$

with

$$u(0, x) = u_0(x), \tag{66b}$$

with measurable $u_0(x) \in [0, 1]$, $x \in R^d$, and $\vec{m} = (m_1, \ldots, m_d) = \sum_z z p(z)$. Its connection with the exclusion processes is given by Theorem 9, which is a paraphrase of Theorems 1.3 and 7.1 in [16]. The first result (70) is the fundamental hydrodynamic limit; the second result (71) is a modified local limit. The original theorems apply to zero range processes as well.

Theorem 9 employs the following notation. We let $\lambda^n$, $n = 1, 2, \ldots$, denote product measures on $Z^d$, with

$$\lambda^n\{\eta : \eta(z) = 1\} = u_{z,n} \quad \text{for } z \in Z^d, \tag{67}$$

where, for all $r > 0$,

$$\int_{|x| < r} |u_{[nx],n} - u_0(x)| \, dx \to 0 \quad \text{as } n \to \infty, \tag{68}$$

and $u_0(x) \in [0, 1]$ is measurable. For a cylinder function $f$ on $\{0, 1\}^{Z^d}$ [i.e., depending only on $\eta(z)$ with $|z| \leq r$, for some $r$], set

$$\hat{f}(\alpha) = E^{\nu_\alpha}[f(\eta)], \quad \alpha \in [0, 1]. \tag{69}$$

Also, set $\tau_z f(\eta) = f(\tau_z \eta)$, where $\tau_z \eta(z') = \eta(z + z')$.

THEOREM 9. *Suppose that the product measures $\lambda^n$, $n = 1, 2, \ldots$, satisfy (67) and (68) for some $u_0(x)$. Then, for any finite open ball $B \subset R^d$, $t \geq 0$ and $\varepsilon > 0$,*

$$P^{\lambda^n}\left\{\left|n^{-d} \sum_{z \in nB} \eta_{nt}(z) - \int_B u(t, x) \, dx\right| > \varepsilon\right\} \to 0 \tag{70}$$

as $n \to \infty$, where $u(t,x)$ is the entropy solution of (66). Moreover, for any cylinder function $f$,

$$(71) \qquad E^{\lambda^n}\left[n^{-d} \sum_{z \in nB} \tau_z f(\eta_{nt})\right] \to \int_B \hat{f}(u(t,x))\,dx$$

as $n \to \infty$.

The limit (71) states that, as $n \to \infty$, $\lambda^n_{nt}$ converges locally to a product measure when viewed on the hydrodynamic scale. With certain monotonicity restrictions on $u_0(x)$, it follows from (71) that, away from discontinuities of $u(t,x)$, $\lambda^n_{nt}$ in fact converges locally to a product measure without any averaging. This is the content of Proposition 3.

Let $C_{v,\theta}$ be the cone with vertex at the origin, pointing in the direction $v \neq 0$, and including all points in $Z^d$ within angle $\theta > 0$ of $v$. We will require that the product measures $\lambda^n$ satisfy

$$(72) \qquad \lambda^n\{\eta : \eta(z) = 1\} \leq \lambda^n\{\eta : \eta(z + z') = 1\}$$

for all $z \in Z^d$ and $z' \in C_{v,\theta}$, for a given choice of $v$ and $\theta$. That is, translation by $C_{v,\theta}$ increases $\lambda^n$ stochastically. After coupling the corresponding processes, application of Theorem 9 implies the following result with a little work.

PROPOSITION 3. *Suppose that the product measures $\lambda^n$, $n = 1, 2, \ldots$, satisfy* (67), (68) *and* (72). *Then,*

$$(73) \qquad \tau_{[nx]}\lambda^n_{nt} \xrightarrow{w} \nu_{u(t,x)} \qquad as\ n \to \infty$$

*for any $t > 0$ and continuity point $x$ of $u(t,\cdot)$, where $u(t,x)$ is the entropy solution of* (66).

We omit a proof of Proposition 3 since it is close to results contained in [10]. Theorem 3 there assumes $u_0(x)$ is continuous rather than requiring (72); a corollary assumes that $\lambda^n = \lambda$ is fixed, with density $\alpha_1$ in one octant in $R^d$ and another density $\alpha_2$ elsewhere. Both results hold for a more general class of particle systems that include the zero range processes as well.

In order to obtain explicit limits in (73), one needs to be able to solve Burgers' equation (66) with its assigned initial data. Fortunately, the equation is degenerate in the sense that its solutions are given by a $(d-1)$-dimensional family of solutions of the one-dimensional Burgers' equation, along lines pointed in the direction $\vec{m} = (m_1, \ldots, m_d)$; this simplifies the computations. To state the result, Proposition 4, we choose an orthonormal basis $v_1, \ldots, v_d$ with $v_1 = \vec{m}/|\vec{m}|$, let $(y, w_2, \ldots, w_d)$ denote the coordinates of $x$ with respect to $v_1, \ldots, v_d$, and set $w = (w_2, \ldots, w_d)$. (If $\vec{m} = 0$, any orthonormal basis can be chosen.) We defer its proof to the Appendix.

PROPOSITION 4. *Suppose $u_0(x)$, $x \in R^d$, is measurable with $u_0(x) \in [0,1]$. Let $u(t,y;w)$ denote the family of entropy solutions of* (63a), *with $m = |\vec{m}|$ and $u_0(y;w) = u_0(x)$, and assume that*

$$u(t,x) \stackrel{def}{=} u(t,y;w) \tag{74}$$

*is jointly measurable in $t$ and $x$. Then, $u(t,x)$ is the entropy solution of* (66a) *with initial data $u_0(x)$.*

We provide here several applications of Propositions 3 and 4, where we explicitly solve for the solutions $u(t,x)$ of (66). In each case, we examine $\lambda_t$ as $t \to \infty$, for initial product measures $\nu_\alpha$ with densities $\alpha(z)$ satisfying

$$\begin{aligned} \alpha(z) &= r \qquad \text{on } A, \\ &= \ell \qquad \text{on } A^c, \end{aligned} \tag{75}$$

for given $\ell, r \in [0,1]$ and $A \subset Z^2$. Analogous examples hold in $Z^d$, $d > 2$. The following examples include two-dimensional analogues of the limits in (62), when viewed away from discontinuities of $u(t,\cdot)$.

EXAMPLE 1. Let $A$ denote the half-space $z^{(1)} \geq cz^{(2)}$ for some $c \geq 0$, where $z = (z^{(1)}, z^{(2)})$, and assume $\vec{m} = (1,0)$. In (67), we can set $u_{z,n} = r$ on $A$ and $u_{z,n} = \ell$ on $A^c$, so that $\lambda^n = \lambda = \nu_\alpha$, for all $n$, where $\alpha$ is given in (75). Since $nA = A$, we choose

$$\begin{aligned} u_0(x) &= r \qquad \text{on } A, \\ &= \ell \qquad \text{on } A^c, \end{aligned} \tag{76}$$

in (68). The cone condition (72) is clearly satisfied. Consequently, by Proposition 3,

$$\tau_{[nx]} \lambda_{nt} \stackrel{w}{\to} \nu_{u(t,x)} \qquad \text{as } n \to \infty, \tag{77}$$

with $u(t,x)$ satisfying (66). Application of Proposition 4 reduces the computation of $u(t,x)$ to solving (63), with $y = x^{(1)} - cx^{(2)}$ and $m = 1$. The one-dimensional problem breaks into two cases, depending on whether or not $\ell < r$. If $\ell < r$, then $u(t,y)$ is given by (64); otherwise, it is given by (65). Note that the resulting solution $u(t,x)$ is measurable, as required for Proposition 4.

EXAMPLE 2. Let $A$ denote the wedge $A = A_1 \cap A_2$, where $A_1$ and $A_2$ are, respectively, the half-spaces where $z^{(1)} \geq cz^{(2)}$ and $z^{(1)} \geq -cz^{(2)}$ for some $c > 0$. As before, we assume that $\vec{m} = (1,0)$. It again follows from Proposition 3 that (77) holds for initial data given by (76). Application of Proposition 4 reduces the problem to solving (63), with $y = x^{(1)} - c|x^{(2)}|$ and $m = 1$. As in Example 1, the problem breaks into two cases depending on whether or not $\ell < r$. Again, $u(t,y)$ is given by either (64) or (65).

EXAMPLE 3. Let $A$ denote the wedge in Example 2, but assume that $\vec{m} = (0, 1)$ here. The limit (77) again holds for the initial data in (76). Application of Proposition 4 reduces the problem to solving (63a) for $m = 1$ and

$$u_0(y) = r \qquad \text{on } A(x^{(1)}),$$
$$= \ell \qquad \text{on } A(x^{(1)})^c,$$

where

$$A(x^{(1)}) = \{y : y \in [-x^{(1)}/c, x^{(1)}/c]\},$$

with $y = x^{(2)}$.

For $x^{(1)} \leq 0$, one trivially obtains $u(t, y) \equiv \ell$. For $x^{(1)} > 0$, we first consider $\ell > r$, with the case $\ell < r$ being analogous. We use standard arguments similar to those on pages 291–303 in [17]. One can check that the behavior of $u(t, y)$ depends on whether (a) $t \leq 2x^{(1)}/c(\ell - r)$ or (b) $t > 2x^{(1)}/c(\ell - r)$. Under (a), $u(t, y)$ is continuous and piecewise linear except at $y = (1 - \ell - r)t + x^{(1)}/c$, with

$$u(t, y) = r \qquad \text{for } y \in [(1 - 2r)t - x^{(1)}/c, (1 - \ell - r)t + x^{(1)}/c],$$
(78)
$$\phantom{u(t, y)} = \ell \qquad \text{for } y \in (-\infty, (1 - 2\ell)t - x^{(1)}/c)$$
$$\cup ((1 - \ell - r)t + x^{(1)}/c, \infty).$$

That is, there is a shock wave emanating from $x^{(1)}/c$ and a rarefaction wave emanating from $-x^{(1)}/c$ at $t = 0$, that first meet at $t_0 = 2x^{(1)}/c(\ell - r)$.

Under (b), $u(t, y)$ is continuous and piecewise linear except at

(79) $$b(t) = (1 - 2\ell)t + 2\sqrt{2(\ell - r)x^{(1)}t/c} - x^{(1)}/c,$$

with

$$u(t, y) = \ell - \sqrt{2(\ell - r)x^{(1)}/ct} \qquad \text{for } y = b(t),$$
(80)
$$\phantom{u(t, y)} = \ell \qquad\qquad\qquad\qquad\quad \text{for } y \in (-\infty, (1 - 2\ell)t - x^{(1)}/c)$$
$$\cup (b(t), \infty).$$

[We are assuming here that $u(t, \cdot)$ is left continuous at $b(t)$, i.e., $b(t)$ lies on the rarefaction wave.] The derivation of (79) and (80) requires some computation. If $b(t)$, $t \geq t_0$, is the position of the discontinuity at which the shock and rarefaction waves meet, then

(81) $$b'(t) = 1 - \ell - u(t, b(t)).$$

This motion is analogous to that of the shock in (64). Moreover, because $(t, b(t))$ lies on the rarefaction wave,

(82) $$u(t, b(t)) = \frac{1}{2t}(t - b(t) - x^{(1)}/c),$$

which is analogous to the stretching of the wave in (65). Substitution of (82) into (81) gives a first-order linear equation having the solution (79); together with (82), this also yields (80). The analogous formulas hold for the case $\ell < r$, with $-y$ being substituted for $y$ in (78) and (80).

In Theorem 4, we showed that if $v \in Z^d \setminus \{0\}$ and $p(\cdot)$ has a finite mean $\vec{m}$, then any $v$-profile measure will satisfy $\langle \vec{m}, v \rangle \geq 0$. Using (70) of Theorem 9, we give a short alternative proof here. As elsewhere in this section, we assume that $p(\cdot)$ has finite range. (The condition $v \in Z^d$ from Section 5 is not needed here, however.)

In order to employ (70), it is useful to be able to omit the assumption that the measures $\lambda^n$ in (67) are product measures. Let the measurable map $U : R_+ \times R^d \to P([0,1])$ (the probability measures on $[0,1]$ equipped with the weak topology) be a (measure-valued) weak limit, under hydrodynamic scaling, of a sequence of exclusion processes with initial measures $\lambda^n$. We will assume any such limit satisfies the regularity condition

$$(83) \qquad \int_{|x|<r} dx \int_{[0,1]} |y - u_0(x)| (U(t,x))(dy) \to 0 \qquad \text{as } t \searrow 0$$

(a.e. in $t$), for all $r > 0$ and a given measurable $u_0(x) \in [0,1]$. Kipnis and Landim [8] comment, on page 199, that one can substitute (83) in Theorem 9 for the assumption that $\lambda^n$ are product measures. (This is only stated for zero range processes restricted to the scaled unit torus, but it is also true in the present setting [11].) We will use this variant of Theorem 9 in the proof of Theorem 10.

THEOREM 10. *Suppose that $\mu$ is a $v$-profile measure that is stationary for the exclusion process. Then, $\langle \vec{m}, v \rangle \geq 0$.*

PROOF. Set $u_{z,n} = \mu\{\eta : \eta(z) = 1\}$ for all $z \in Z^d$ and $n \in Z^+$. Because $\mu$ is a $v$-profile measure, $u_{z,n}$ satisfies (68) with

$$(84) \qquad \begin{aligned} u_0(x) &= 1 \qquad \text{for } \langle x, v \rangle \geq 0, \\ &= 0 \qquad \text{for } \langle x, v \rangle < 0. \end{aligned}$$

Likewise, since $\mu_t = \mu$, for all $t$, (83) trivially holds for limits $U(t,x)$ of the exclusion process under hydrodynamic scaling. So, (70) is satisfied for the entropy solution $u(t,x)$ of (66). Again, since $\mu_t = \mu$, one must have, for all $t$,

$$(85) \qquad \begin{aligned} u(t,x) &= 1 \qquad \text{for } \langle x, v \rangle \geq 0, \\ &= 0 \qquad \text{for } \langle x, v \rangle < 0. \end{aligned}$$

This implies $\langle \vec{m}, v \rangle \geq 0$. Otherwise, as in Proposition 4, let $u(t, y; w)$, $w \in R^{d-1}$, denote the entropy solutions of (63a), with $m = |\vec{m}|$. If $\langle \vec{m}, v \rangle < 0$, then (63b) is satisfied with $\ell = 1$ and $r = 0$, for each $w$. One obtains the solution (65), which contradicts (85). $\square$

We note that one can also demonstrate Theorem 10 by using the constant sequence of product initial measures with $u_{z,n} = \mathbb{1}_{\langle z,v \rangle \geq 0}$, rather than appealing to (83). The argument involves truncations, and so is not as quick as before.

**9. An example with local behavior depending on $p(\cdot)$.** Results such as Theorem 9 and Proposition 3, that use hydrodynamic limits, provide substantial insight into the asymptotic behavior of exclusion processes on $Z^d$. A more refined analysis is needed on the original scale, however, to obtain information about the exclusion process at the shocks of solutions of the corresponding Burgers' equation. For instance, the limit of $\lambda_t$ in (62), for the case where $0 < \ell < r$ and $\ell + r = 1$, is given by a mixture of the product measures $\nu_\ell$ and $\nu_r$, since the shock can randomly be on either side of 0.

We give here another type of example, that illustrates how the asymptotic behavior can depend on $p(\cdot)$ itself, and not just on its mean $\vec{m}$, if $d \geq 2$. For this, we examine exclusion processes on $Z^2$ with initial measure

$$\begin{aligned}
\lambda\{\eta : \eta(z) = 1\} &= 1 \quad \text{for } z \in A, \\
&= 0 \quad \text{for } z \in A^c,
\end{aligned} \tag{86}$$

where $A = \{z = (z^{(1)}, z^{(2)}) : z^{(1)} \geq 0 \text{ and } z^{(2)} \geq 0\}$, and with random walk kernels

$$\begin{aligned}
p(z) &= p_1 &&\text{for } z = e_1, e_2, \\
&= q_1 = \tfrac{1}{2} - p_1 &&\text{for } z = -e_1, -e_2,
\end{aligned} \tag{87}$$

for $p_1 > 1/4$, and

$$\begin{aligned}
p(z) &= p_2 &&\text{for } z = e_1, e_2, \\
&= q_2 = 1 - 2p_2 &&\text{for } z = 0.
\end{aligned} \tag{88}$$

Under (87), $\vec{m} = (2p_1 - \tfrac{1}{2}, 2p_1 - \tfrac{1}{2})$, and under (88), $\vec{m} = (p_2, p_2)$, which are, of course, equal for $p_2 = 2p_1 - \tfrac{1}{2}$. Since particles can never move under (86) and (88), $\lambda_t = \lambda$ for all $t$, in this case. As Theorem 11 shows, the behavior under (86) and (87) is quite different.

THEOREM 11. *Assume that an exclusion process has transition probabilities $p(\cdot)$ satisfying (87) and initial measure satisfying (86). Then, for each $j$,*

$$P^\lambda\{\eta : \eta_t(z) = 1 \text{ for some } z^{(1)} + z^{(2)} \leq j\} \to 0 \quad \text{as } t \to \infty. \tag{89}$$

A little thought shows what, in principle, must be occurring for the process in Theorem 11. Particles randomly move into vertical and horizontal strips just to the left and below the $z^{(2)}$ and $z^{(1)}$ axes, along which they can drift up and to the right. The local departure of particles continues over time, emptying the plane locally and producing (89).

The demonstration of Theorem 11 takes up the remainder of the section. Much of the work is contained in Lemma 7, which employs the following notation. We define $\lambda^j$ as in (86), but where $A$ is replaced by $D_j = \{z : z^{(1)} + z^{(2)} > j\}$, and set

$$\Delta_j = \{z : z^{(1)} \geq -j,\ z^{(2)} \geq -j,\ z^{(1)} + z^{(2)} \leq j\}$$

and

$$V_j = D_j^c \cap \{z : z^{(1)} < -j \text{ or } z^{(2)} < -j\}.$$

In words, $\Delta_j$ is the triangle with vertices at $(-j, -j)$, $(-j, 2j)$ and $(2j, -j)$, and $V_j$ is the half-plane below $z^{(1)} + z^{(2)} = j$, excluding $\Delta_j$. Inequalities (90) and (91) in Lemma 7 give upper bounds on the probabilities of there being particles in $V_j$ and $\Delta_j$ for the exclusion process starting at $\lambda$ and $\lambda^{2j}$. From now on, we consider $p_1 > 1/4$ in (87) to be fixed.

LEMMA 7. *For appropriate $c > 0$,*

(90) $$P^\lambda\{\eta : \eta_t(z) = 1 \text{ for some } z \in V_j\} \leq c(q_1/p_1)^{j/3}$$

*and*

(91) $$P^{\lambda^{2j}}\{\eta : \eta_t(z) = 1 \text{ for some } z \in \Delta_j\} \leq c(q_1/p_1)^{j/3}$$

*for all $t$.*

PROOF. In order to demonstrate (90), it suffices to show the inequality with $V_j$ replaced by $V_{j,1}$ and by $V_{j,2}$, where

$$V_{j,1} = V_j \cap \{z : z^{(1)} \leq z^{(2)}\}, \qquad V_{j,2} = V_j \cap \{z : z^{(1)} \geq z^{(2)}\}.$$

By symmetry, it suffices to do this for just $V_{j,1}$.

Our main estimate will be a strengthened version of (44) of Lemma 3, with an explicit rate of decay. For this, we employ the notation

$$L(z) = \{z' \in Z^2 : (z')^{(2)} = z^{(2)}\},$$

$$\ell_{\lambda'}^-(z) = E^{\lambda'}[\#z' : (z')^{(1)} \leq z^{(1)},\ (z')^{(2)} = z^{(2)} \text{ and } \eta(z') = 1]$$

and

$$\ell_{\lambda'}^+(z) = E^{\lambda'}[\#z' : (z')^{(1)} \geq z^{(1)},\ (z')^{(2)} = z^{(2)} \text{ and } \eta(z') = 0],$$

where $\lambda'$ is defined as in (86), but with $A$ replaced by $A' = \{z : z^{(1)} \geq 0\}$. The terms $L$, $\ell_{\lambda'}^-$ and $\ell_{\lambda'}^+$ are special cases of those in Section 7, with $d = 2$ and $v = (1, 0)$ here. We let $\nu = \nu_\alpha$ denote the product measure with

$$\pi(z) = (p_1/q_1)^{z^{(1)}} \tag{92}$$

[where $\pi(z) = \alpha(z)/(1 - \alpha(z))$].

Using Theorem 2, it is easy to check that $\nu$ is stationary. Clearly, $\nu$ and $\lambda'$ are $v$-homogeneous measures for $v = (1, 0)$. They also satisfy both parts of (42). Setting $\mu = \nu$ and $\lambda = \lambda'$, the assumptions of Lemma 3 are therefore satisfied. Consequently, by (44),

$$\ell_{\lambda'_t}^-(z) \to 0 \qquad \text{as } z^{(1)} \to -\infty$$

uniformly in $t$. Using the explicit form of $\pi(\cdot)$ in (92), this can be strengthened by setting $\varepsilon = c_1(q_1/p_1)^n$, with appropriate $c_1 > 0$, in the first two displays of the demonstration of (44). Combining these two inequalities, as before, implies that

$$\ell_{\lambda'_t}^-(z) < 2c_1(q_1/p_1)^n \tag{93}$$

for $z^{(1)} = -2n$ and all $t$. Since $A \subset A'$, it follows that

$$\ell_{\lambda_t}^-(z) < 2c_1(q_1/p_1)^n \tag{94}$$

for $z^{(1)} = -2n$ and all $t$.

The region $V_{j,1}$ is the union of the horizontal line segments $(-\infty, r(z^{(2)})]$, $z^{(2)} \in Z$, with

$$r(z^{(2)}) = -j - 1 \qquad \text{for } z^{(2)} \in [-j, j],$$
$$= -|z^{(2)}| \qquad \text{for } z^{(2)} \notin [-j, j].$$

Together with (94), this implies

$$E^\lambda[\#z \in V_{j,1} : \eta(z) = 1] \leq (2j+1)2c_1(q_1/p_1)^{[j/2]} + 2 \sum_{n=j+1}^\infty 2c_1(q_1/p_1)^{[n/2]}$$

$$\leq c(q_1/p_1)^{j/3}$$

for large enough $c$. This demonstrates (90).

The reasoning for (91) is similar. Here, we set

$$L(z) = \{z' \in Z^2 : (z')^{(1)} - (z')^{(2)} = z^{(1)} - z^{(2)}\}. \tag{95}$$

The expectations $\ell_{\lambda^j}^-(z)$ and $\ell_{\lambda^j}^+(z)$ are defined as before, but with the inner equality in (95) replacing $(z')^{(2)} = z^{(2)}$. The vector $v = (1, 1)$ replaces $(1, 0)$, and we let $\nu^j = \nu_\alpha$ denote the product measure with

$$\pi(z) = (p_1/q_1)^{z^{(1)} + z^{(2)} - j}. \tag{96}$$

The measures $\nu^j$ are stationary; $\nu^j$ and $\lambda^j$ are $v$-homogeneous and satisfy both parts of (42). Setting $\mu = \nu^{2j}$ and $\lambda = \lambda^{2j}$, it follows as before that, for appropriate $c_1 > 0$,

$$\ell^-_{\lambda^{2j}_t}(z) < 2c_1(q_1/p_1)^n$$

for $z^{(1)} + z^{(2)} = 2j - 2n$ and all $t$. The region $\Delta_j$ is the union of $2j + 1$ line segments pointing in the direction $v$ with right endpoint $z^{(1)} + z^{(2)} = j$ in each case. Setting $n = [j/2]$ in each case implies that

$$E^{\lambda^{2j}}[\#z \in \Delta_j : \eta(z) = 1] \leq (2j+1)2c_1(q_1/p_1)^{[j/2]} \leq c(q_1/p_1)^{j/3}$$

for large enough $c$. This demonstrates (91).  □

The proof of Theorem 11 applies Lemma 7, together with a coupling argument that employs the translation invariance of $D_{2j}$ under $(-1, 1)$. The basic idea is that since

$$\{z : z^{(1)} + z^{(2)} \leq j\} = V_j \cup \Delta_j,$$

(89) will follow from (90) and (91), if the initial measure $\lambda^{2j}$ in (91) can be replaced by $\lambda$ for large $t$. The difference $A - D_{2j}$ of the sets of sites initially occupied by the corresponding processes is finite. So, it is reasonable to expect that the particles originally there will eventually dissipate away from the finite region $\Delta_j$ as $t \to \infty$, perhaps by moving more or less in the direction $(-1, 1)$.

The coupling we will use, like that after (50), is an extension of the coupling $\xi_t = (\eta_t, \zeta_t)$ to a coupling $\vec{\xi}_t = (\vec{\eta}_t, \zeta_t)$ with multiple coordinates. Here, $\vec{\eta} = (\eta^0, \ldots, \eta^{N-1})$, with $N = [19j^4/\varepsilon]$, for a given $\varepsilon > 0$. The coordinate processes are assumed to have deterministic initial measures $\mu^{2j,n}$, $n = 0, 1, \ldots, N-1$, and $\lambda^{2j}$. The measures $\lambda^{2j}$ are defined as before; $\mu^{2j,n}$ are defined as in (86), but with $A$ replaced by $D_{2j,n} = D_{2j} \cup \Delta'_{2j,n}$, where

$$\Delta'_{j,n} = \{z : z^{(1)} \geq 0,\ z^{(2)} \geq 0,\ z^{(1)} + z^{(2)} \leq j\} + (-2jn, 2jn)$$

is the translation, by $(-2jn, 2jn)$, of the triangle with vertices at $(0,0)$, $(0, j)$ and $(j, 0)$. Let $\vec{\gamma}$ denote the joint initial measure of the coupling.

We note that $\mu^{2j,n}$ is the translate of $\mu^{2j,0}$ by $z = (-4jn, 4jn)$, and so $\mu^{2j,n}_t$ is also the translate of $\mu^{2j,0}_t$ by $z$; since $A \subset D_{2j,0}$, $\mu^{2j,0}_t$ stochastically dominates $\lambda$. Also, since $D_{2j} \subset D_{2j,n}$, for all $j$ and $n$, the discrepancies for a pair $(\eta^n_t, \zeta_t)$ are always of type $(1, 0)$. Their number remains constant over time, which is $2j^2 + j$, the number of sites in $\Delta'_{j,n}$. The distribution of such discrepancies, at a time $t$, is the same as for $n = 0$, up to a translation by $z$. We also set

$$\Delta_{j,n} = \Delta_j + (-4jn, 4jn),$$

and note that $\Delta_j^{n_1} \cap \Delta_j^{n_2} = \varnothing$ for $n_1 \neq n_2$.

PROOF OF THEOREM 11. Set
$$F_t^n = \{\vec{\xi} : \eta_t^n(z) = 1,\ \zeta_t(z) = 0 \text{ for some } z \in \Delta_{j,n}\},$$

for $n = 0, 1, \ldots, N - 1$. We will show that, for fixed $\varepsilon > 0$,

(97) $$\sum_{n=0}^{N-1} P^{\vec{\gamma}}(F_t^n) \leq \varepsilon N$$

for large enough $t$. It follows from the above comments on discrepancies that the probabilities of the sets in (97) are all the same, and so $P^{\vec{\gamma}}(F_t^0) \leq \varepsilon$. Consequently,

$$P^{\mu^{2j,0}}\{\eta : \eta_t(z) = 1 \text{ for some } z \in \Delta_j\} - P^{\lambda^{2j}}\{\eta : \eta_t(z) = 1 \text{ for some } z \in \Delta_j\} \leq \varepsilon,$$

for such $t$. Choosing $\varepsilon \leq c(q_1/p_1)^{j/3}$, it follows from this, (91) and $A \subset D_{2j,0}$, that, for given $j$ and large enough $t$,

$$P^\lambda\{\eta : \eta_t(z) = 1 \text{ for some } z \in \Delta_j\} \leq 2c(q_1/p_1)^{j/3}.$$

Along with (90), this implies that

$$P^\lambda\{\eta : \eta_t(z) = 1 \text{ for some } z^{(1)} + z^{(2)} \leq j\} \leq 3c(q_1/p_1)^{j/3},$$

which implies (89) and, hence, Theorem 11.

In order to demonstrate (97), we set

$$X_t(\vec{\xi}) = \#\{n : \vec{\xi} \in F_t^n\}.$$

We claim that, on $X_t > 9j^4$, there must exist $n_1, n_2$ and $z_1, z_2$, with $z_i \in \Delta_{j,n_i}$, so that

(98) $$\eta_t^{n_1}(z_1) = \eta_t^{n_2}(z_2) = 1 \quad \text{and} \quad \eta_t^{n_1}(z_2) = \eta_t^{n_2}(z_1) = 0.$$

To show this, let $n'_1, \ldots, n'_{3j^2}$ denote the first $3j^2$ indices for which $\vec{\xi} \in F_t^{n'_i}$. Then,

$$\left| \bigcup_{i=1}^{3j^2} \{z : \eta_t^{n'_i}(z) = 1,\ \zeta_t(z) = 0\} \right| \leq 3j^2 |\Delta_j| \leq 9j^4.$$

So, there is an index $n_1$, with

(99) $$\eta_t^{n_1}(z_1) = 1 \quad \text{and} \quad \eta_t^{n'_i}(z_1) = \zeta_t(z_1) = 0$$

for some $z_1 \in \Delta_{j,n_1}$ and all $i$. For some $n_2 \in \{n'_1, \ldots, n'_{3j^2}\}$,

$$\{z \in \Delta_{j,n_2} : \eta_t^{n_1}(z) = 1 \text{ and } \zeta_t(z) = 0\} = \varnothing,$$

again since $|\Delta_j| \leq 3j^2$. One can therefore choose $z_2 \in \Delta_{j,n_2}$ so that

$$\eta_t^{n_2}(z_2) = 1 \quad \text{and} \quad \eta_t^{n_1}(z_2) = \zeta_t(z_2) = 0.$$

Together with (99), this implies (98).

When the event in (98) occurs, the pair $(\eta_t^{n_1}, \eta_t^{n_2})$ has discrepancies of opposite types at $z_1$ and $z_2$. Since $z_i \in \Delta_{j,n_i}$, one has $|z_2 - z_1| \leq 4\sqrt{2}jN$. In the same manner as in the proof of Lemma 4, one can employ Lemma 3.1 of [1] to obtain a uniform lower bound $\delta > 0$ on the probability of these or other discrepancies of opposite type meeting over $(t, t+1]$. Also, there are not more than $3j^2N$ sites at which $\vec{\xi}_0$ has discrepancies. Therefore, the event $X_t > 9j^4$ cannot occur with probability at least $\varepsilon/2$ at more than $6j^2N/\delta\varepsilon$ times spaced at least distance 1 apart. It follows that, for given $\varepsilon > 0$ and sufficiently large $T_\varepsilon$,

$$P^{\vec{\gamma}}\{X_t > 9j^4\} < \varepsilon/2 \tag{100}$$

for $t \geq T_\varepsilon$.

One has

$$\sum_{n=0}^{N-1} P^{\vec{\gamma}}(F_t^n) = E^{\vec{\gamma}}[X_t] = E^{\vec{\gamma}}[X_t; X_t > 9j^4] + E^{\vec{\gamma}}[X_t; X_t \leq 9j^4].$$

Since $X_t \leq N$, this is, by (100) and our choice of $N$,

$$\leq \tfrac{1}{2}\varepsilon N + 9j^4 \leq \varepsilon N,$$

for small $\varepsilon$. This implies (97), and completes the proof of Theorem 11. □

**10. Open problems.** An important open problem in the context of this paper is to determine completely the set of extremal stationary measures for an irreducible translation invariant exclusion process on $Z^d$, $d > 1$, with $\sum_x |x|p(x) < \infty$. A complete characterization is out of reach at the present time. It is therefore useful to isolate particular parts of this problem that are deserving of attention. Here are some (related) statements that we think are probably true, but cannot prove. [Some may require higher moment assumptions on $p(\cdot)$.]

1. All extremal inhomogeneous stationary measures are $v$-homogeneous for some $v \in R^d$ satisfying $\langle \sum_x xp(x), v \rangle > 0$.
2. All extremal inhomogeneous stationary measures are $v$-profile for some $v \in R^d$ satisfying $\langle \sum_x xp(x), v \rangle > 0$. (We are less sure about this statement.)
3. For each $v \in R^d$ satisfying $\langle \sum_x xp(x), v \rangle > 0$, there is a continuous one parameter family of inhomogeneous extremal stationary $v$-profile measures.

4. If $\sum_x x p(x) = 0$, then the extremal stationary measures are exactly $\{\nu_\rho, 0 \leq \rho \leq 1\}$. This is an old problem—it appears as open problem #5 on page 416 in [13].
5. For every choice of $p(\cdot)$ with nonzero mean vector, there exists an extremal stationary measure that is not a product measure. [We do not know that this is the case for *any* $p(\cdot)$ if $d \geq 2$.]

In one dimension, #1 above is trivial, #2 is true, #3 is true with extra assumptions, except that the parameter is discrete rather than continuous, #4 is true, and #5 is true under extra assumptions. (For the extra assumptions required, see the discussion of the one-dimensional system in Section 1.)

Problems #1, #2 and #3 are interesting and open if $Z^d$ is replaced by a "ladder" of the form $Z \times \{1, \ldots, N\}$. A possible analogue of the $v$-homogeneous and $v$-profile properties in this context is the following: One can view $\{1, \ldots, N\}$ as a cycle, and say that a measure is rotationally homogeneous if it is invariant under "rotations" of $Z \times \{1, \ldots, N\}$, and is profile if it is rotationally homogeneous and has density tending to 1 in one direction and to 0 in the other. Problem #4 can be handled in the case of a ladder by a small modification of the one-dimensional proof. It is likely that #5 for the ladder can be done using the approach in [2].

It may be useful to recall the known characterization of the class of all extremal stationary measures for continuous time independent particle systems on $Z^d$ with transition kernel $p(\cdot)$, since this might shed some light on what to expect for the exclusion process. By Theorem 4.12 in [12], the extremal stationary measures are exactly the Poisson fields on $\{0, 1, 2, \ldots\}^{Z^d}$ in which the number of particles at $x$ has mean $m(x)$, where $m(x) \geq 0$ for all $x$ and $\sum_x m(x) p(y - x) = m(y)$ for all $y$. By the Choquet–Deny theorem, the extremal functions $m$ satisfying these properties are the pure exponentials

$$m(x) = e^{\langle x, v \rangle},$$

for $v$ such that $\sum_x e^{-\langle x, v \rangle} p(x) = 1$. Nonconstant $m(\cdot)$ of this type can only exist if $\sum_x x p(x) \neq 0$, and then the corresponding $v$ must satisfy $\langle \sum_x x p(x), v \rangle > 0$. Of course, nonextremal $m$'s can correspond to extremal stationary Poisson fields.

Note that the Poisson fields corresponding to the pure exponentials are $v$-homogeneous (and even satisfy a property analogous to being $v$-profile), while those corresponding to sums of exponentials are not. In one dimension, there are no extremal stationary measures for the exclusion process corresponding to sums of exponentials—they are all either homogeneous or $v$-profile. Perhaps this is true in higher dimensions as well. This has a bearing on problem #1 above.

In the context of the convergence results of Section 7, it is not clear to us what assumptions on $p(\cdot)$ and $v$ are, in general, necessary in order for the conditions (42) and (48) of Theorem 8 to hold. It seems reasonable to guess that the existence of a one-parameter family $\mu^\alpha$ of extremal stationary $v$-profile measures, as in (48), is in some sense generic.

The tail condition (42) should presumably hold for the $v$-profile measures of a large class of $p(\cdot)$. A related question arose in [1] in $d=1$. As pointed out there, presently there are no known examples of extremal stationary profile measures that are not blocking measures. Our condition (42) corresponds to (1.3) in that paper, which immediately implies that a measure is a blocking measure. It is also not known what conditions on $p(\cdot)$ are necessary for all stationary blocking measures to satisfy (1.3); it is suggested there that a third moment assumption on $p(\cdot)$ might be the correct condition.

## APPENDIX

In Section 8, we examined the entropy solutions of the $d$-dimensional Burgers' equation in (66). The equation is a specific case of the scalar conservation law

$$\text{(A1a)} \qquad \frac{\partial u}{\partial t} + \sum_{i=1}^{d} m_i \frac{\partial}{\partial x_i} f_i(u) = 0,$$

with $f_i \in C^2(R)$, and having bounded initial data

$$\text{(A1b)} \qquad u(0,x) = u_0(x), \qquad x \in R^d.$$

We briefly review here basic existence, uniqueness and stability results for the entropy solutions of (A1). We then give the quick proof of Proposition 4. More detail can be found on conservation laws in one dimension in Chapters 15 and 16 in [17] and Chapter 3 in [4], and on conservation laws in $d$-dimensions in [9] and in Appendix 2 in [8].

Classical solutions to (A1) need not exist for all $t$, even when the initial data $u_0(x)$ is smooth. One therefore typically works with weak solutions of (A1). That is, for every function $g: R_+ \times R^d \to R$ of class $C_K^{1,1}(R_+ \times R^d)$, $u: R_+ \times R^d \to R$ is required to satisfy

$$\text{(A2)} \qquad \int_0^\infty \int_{R^d} \left( u(t,x) \frac{\partial}{\partial t} g(t,x) + \sum_{i=1}^{d} m_i f_i(u(t,x)) \frac{\partial}{\partial x_i} g(t,x) \right) dx\, dt$$
$$+ \int_{R^d} u_0(x) g(0,x)\, dx = 0,$$

where $f_i$ and $u_0$ are as in (A1). ($C_K^{1,1}$ denotes the continuous functions with compact support and one continuous derivative in both time and space.)

Solutions $u(t,x)$ of (A2) can be modified on sets of measure 0 on $R_+ \times R^d$ without affecting the left-hand side of (A2), and so are certainly not pointwise unique. Because of the possibility of shocks, they also need not be a.e. unique. (See, e.g., Chapter 15 in [17] for examples.) The physically meaningful solution of (A2) is given by its (a.e.) unique entropy solution. The following entropy condition is due to Kruzkov [9], and states that, for every positive function $g$ of class $C_K^{1,1}(R_+ \times R^d)$ and every $a > 0$, $u(t,x)$ satisfies

(A3a)
$$\int_0^\infty \int_{R^d} \left\{ |u(t,x) - a| \frac{\partial}{\partial t} g(t,x) \right. \\ \left. + \sum_{i=1}^d m_i |f_i(u(t,x)) - f_i(a)| \frac{\partial}{\partial x_i} g(t,x) \right\} dx\, dt \geq 0,$$

and for all $r > 0$,

(A3b)
$$\int_{|x|<r} |u(t,x) - u_0(x)|\, dx \to 0 \qquad \text{as } t \searrow 0$$

(a.e. in $t$), where $f_i$ and $u_0$ are chosen as in (A1). A solution $u(t,x)$ of (A2) satisfying (A3) is said to be an *entropy* solution. (Kruzkov [9] includes a more general family of equations in place of (A2).)

Kruzkov [9] proved the existence and uniqueness of entropy solutions:

THEOREM A1. *For every bounded measurable $u_0(x)$, there exists a unique entropy solution of* (A2).

The above uniqueness is a consequence of an $L^1$ stability result in [9] that compares entropy solutions with different initial data.

The theory for solutions of (A2) in one dimension preceded that for general $d$ (see, e.g., [15]). When, in addition, $f'' < 0$, there is a version of the entropy solution satisfying

(A4) $$u(t,y) - u(t, y+a) \leq Ca/t,$$

for $a, t > 0$, where $C > 0$ depends only on $f$ and $\|u_0(y)\|_\infty$ (see, e.g., page 266 in [17]). In particular, only increasing jumps can occur as one moves from left to right, and for a given $t$, $u(t, \cdot)$ is locally of bounded variation. In this setting, one can employ (A4) as the entropy condition in place of (A3).

We conclude with the proof of Proposition 4 of Section 8. The argument is a straightforward application of Fubini's theorem.

PROOF OF PROPOSITION 4. We need to show that (A2) and (A3) hold for $f_i(u) = f(u) = u(1-u)$ and all $g(t,x)$, if $u(t,x)$ satisfies (74). Using the coordinates of the orthonormal basis introduced before Proposition 4 and Fubini's theorem, we can rewrite the left side of (A2) as

$$\int_{R^{d-1}} \left\{ \int_0^\infty \int_R \left( u(t,y;w) \frac{\partial}{\partial t} g(t,y;w) \right. \right.$$
$$\left. + mf(u(t,y;w)) \frac{\partial}{\partial y} g(t,y;w) \right) dy\, dt \tag{A5}$$
$$\left. + \int_{R^{d-1}} \int_R u_0(y;w) g(0,y;w)\, dy \right\} dw,$$

where $g(t,y;w) = g(t,x)$. Since $u(t,y;w)$ is an entropy solution of (63a) and $g(\cdot,\cdot;w) \in C_K^{1,1}(R_+ \times R)$ for each $w$, the quantity inside the braces always equals 0 and, hence so does the entire integral. So, $u(t,x)$ satisfies (A2). The argument for (A3) is the same. □

We note that the joint measurability of $u(t,x)$ in $t$ and $x$ was needed in order for (A2) to make sense and to apply Fubini's theorem. By altering each one-dimensional solution $u(t,y;w)$ at a single point $(t_w, y_w)$, one can produce nonmeasurable $u(t,x)$. So, one cannot arbitrarily choose the version of $u(t,y;w)$ for each $w$. We also point out that the above proof does not hold when $f_i$ depends on $i$, since one would need to factor the functions outside of the summation in order to obtain the analogue of (A5).

**Acknowledgment.** The authors thank C. Landim for his comments on hydrodynamic limits.

## REFERENCES

[1] BRAMSON, M., LIGGETT, T. M. and MOUNTFORD, T. (2002). Characterization of stationary measures for one-dimensional exclusion processes. *Ann. Probab.* **30** 1539–1575. MR1944000
[2] BRAMSON, M. and MOUNTFORD, T. (2002). Stationary blocking measures for one-dimensional nonzero mean exclusion processes. *Ann. Probab.* **30** 1082–1130. MR1920102
[3] CHOQUET, G. and DENY, J. (1960). Sur l'equation de convolution $\mu = \mu * \sigma$. *C. R. Acad. Sci. Paris* **250** 799–801. MR0119041
[4] EVANS, L. (1998). *Partial Differential Equations.* Amer. Math. Soc., Providence, RI. MR1625845
[5] FERRARI, P. A., LEBOWITZ, J. L. and SPEER, E. (2001). Blocking measures for asymmetric exclusion processes via coupling. *Bernoulli* **7** 935–950. MR1873836
[6] GEORGII, H. O. (1979). *Canonical Gibbs Measures. Lecture Notes in Math.* **760**. Springer, Berlin. MR0551621
[7] JUNG, P. (2003). Extremal reversible measures for the exclusion process. *J. Statist. Phys.* **112** 165–191. MR1991035
[8] KIPNIS, C. and LANDIM, C. (1999). *Scaling Limits of Interacting Particle Systems.* Springer, Berlin. MR1707314

 [9] KRUZKOV, S. N. (1970). First order quasilinear equations in several independent variables. *Mat. Sb.* **123** 228–255. (English translation in *Math. USSR Sb.* **10** 217–243.) MR0267257
[10] LANDIM, C. (1993). Conservation of local equilibrium for attractive particle systems on $Z^d$. *Ann. Probab.* **21** 1782–1808. MR1245290
[11] LANDIM, C. (2003). Personal communication.
[12] LIGGETT, T. M. (1978). Random invariant measures for Markov chains and independent particle systems. *Z. Wahrsch. Verw. Gebiete* **45** 297–313. MR0511776
[13] LIGGETT, T. M. (1985). *Interacting Particle Systems*. Springer, New York. MR0776231
[14] LIGGETT, T. M. (1999). *Stochastic Interacting Systems*: *Contact, Voter and Exclusion Processes*. Springer, Berlin. MR1717346
[15] OLEINIK, O. A. (1957). Discontinuous solutions of nonlinear differential equations. *Uspekhi Mat. Nauk* (N.S.) **12** 3–73. [English translation in *Amer. Math. Soc. Transl. Ser.* (*2*) **26** (1963) 95–172.] MR0094541
[16] REZAKHANLOU, F. (1991). Hydrodynamic limit for attractive particle systems on $Z^d$. *Comm. Math. Phys.* **140** 417–448. MR1130693
[17] SMOLLER, J. (1983). *Shock Waves and Reaction–Diffusion Equations*. Springer, New York. MR0688146
[18] SPITZER, F. (1970). Interaction of Markov processes. *Adv. Math.* **5** 246–290. MR0268959

SCHOOL OF MATHEMATICS
UNIVERSITY OF MINNESOTA
206 CHURCH ST. SE
MINNEAPOLIS, MINNESOTA 55455
USA
E-MAIL: bramson@math.umn.edu
URL: www.math.umn.edu/pacim/bramson.htm

DEPARTMENT OF MATHEMATICS
UNIVERSITY OF CALIFORNIA,
 LOS ANGELES
405 HILGARD AVE.
LOS ANGELES, CALIFORNIA 90095
USA
E-MAIL: tml@math.ucla.edu
URL: www.math.ucla.edu/~tml/